\documentclass{article}
\RequirePackage[leqno]{amsmath}[2.14]

\usepackage{amssymb}
\usepackage{mathtools}
\usepackage{tabularx}
\usepackage{picture,slashbox}
\usepackage{graphicx}
\usepackage{epsfig,epsf,psfrag,epstopdf,subfigure}
\usepackage{color}
\usepackage{xcolor}
\usepackage{comment}
\usepackage{algorithm}
\usepackage{algorithmic}
\usepackage{url}
\usepackage{enumitem}

\usepackage{mcode}

\usepackage{caption} 
\usepackage{array,booktabs}
\usepackage{siunitx}
\RequirePackage{xr-hyper}[6.00]
\RequirePackage{ifpdf}[2.3]
\RequirePackage{hyperref}[6.83]
\RequirePackage{xcolor}[2.11]
\colorlet{inlinkcolor}{green!50!black}
\colorlet{exlinkcolor}{red!50!black}
\hypersetup{
  colorlinks = true,
  allcolors = inlinkcolor,
  urlcolor = exlinkcolor,
}

\RequirePackage[capitalize,nameinlink]{cleveref}[0.19]

\crefname{section}{section}{sections}
\crefname{subsection}{subsection}{subsections}
\Crefname{section}{Section}{Sections}
\Crefname{subsection}{Subsection}{Subsections}

\Crefname{figure}{Figure}{Figures}

\crefformat{equation}{\textup{#2(#1)#3}}
\crefrangeformat{equation}{\textup{#3(#1)#4--#5(#2)#6}}
\crefmultiformat{equation}{\textup{#2(#1)#3}}{ and \textup{#2(#1)#3}}
{, \textup{#2(#1)#3}}{, and \textup{#2(#1)#3}}
\crefrangemultiformat{equation}{\textup{#3(#1)#4--#5(#2)#6}}%
{ and \textup{#3(#1)#4--#5(#2)#6}}{, \textup{#3(#1)#4--#5(#2)#6}}{, and \textup{#3(#1)#4--#5(#2)#6}}

\Crefformat{equation}{#2Equation~\textup{(#1)}#3}
\Crefrangeformat{equation}{Equations~\textup{#3(#1)#4--#5(#2)#6}}
\Crefmultiformat{equation}{Equations~\textup{#2(#1)#3}}{ and \textup{#2(#1)#3}}
{, \textup{#2(#1)#3}}{, and \textup{#2(#1)#3}}
\Crefrangemultiformat{equation}{Equations~\textup{#3(#1)#4--#5(#2)#6}}%
{ and \textup{#3(#1)#4--#5(#2)#6}}{, \textup{#3(#1)#4--#5(#2)#6}}{, and \textup{#3(#1)#4--#5(#2)#6}}

\crefdefaultlabelformat{#2\textup{#1}#3}

\numberwithin{equation}{section}
\numberwithin{table}{section}
\numberwithin{figure}{section}

\newcommand{\be}{\begin{equation}}
\newcommand{\ee}{\end{equation}}
\newcommand{\ba}{\begin{aligned}}
\newcommand{\ea}{\end{aligned}}
\newcommand{\bea}{\begin{eqnarray}}
\newcommand{\eea}{\end{eqnarray}}

\title{A fast Gauss transform in one dimension using sum-of-exponentials approximations}

\author{Shidong Jiang\thanks{Department of Mathematical Sciences,
    New Jersey Institute of Technology, Newark, New Jersey 07102
    ({\tt shidong.jiang@njit.edu}).}
}

\begin{document}
\maketitle

\begin{abstract}
  We present a fast Gauss transform in one dimension using nearly optimal
  sum-of-exponentials approximations of the Gaussian kernel.
  For up to about ten-digit accuracy, the approximations are obtained via best
  rational approximations of the exponential function on the negative real axis.
  As compared with existing fast Gauss transforms, the algorithm is straightforward for
  parallelization and very simple to implement, with only twenty-four lines
  of code in {\sc MATLAB}. The most expensive part of the algorithm is on the evaluation
  of complex exponentials, leading to three to six complex exponentials FLOPs per point
  depending on the desired precision. The performance of the algorithm is illustrated
  via several numerical examples.
\end{abstract}



%%%%%%%%%%%%%%%%%%%%%%%%%%%%%%%%%%%%%%%%%%%%%%
\section{Introduction} \label{introduction}
%%%%%%%%%%%%%%%%%%%%%%%%%%%%%%%%%%%%%%%%%%%%%%
In this paper, we consider the evaluation
of the Gauss transform in one dimension:
\be
u_i=\sum_{j=1}^N G(x_i-y_j;\delta)\alpha_j,
\qquad x_i, y_j, \alpha_j \in \mathbb{R}, \qquad i=1,\ldots M,
\label{gt}
\ee
where the Gaussian kernel is given by the formula
\be
G(x;\delta)=e^{-\frac{x^2}{4\delta}}.
\label{gkernel}
\ee
Since the Gaussian kernel is the density function of the normal
distribution, the Gauss transform has broad applications
in all areas of statistics such as the kernel method \cite{elgammal2003},
Gaussian processes \cite{waaij2016}, nonparametric statistic models \cite{geman1982},
density estimation \cite{silverman1986}, etc.
Up to a prefactor, the Gaussian kernel is the fundamental
solution of the heat equation. Thus, the Gauss transform
also appears ubiquitously in integral equation methods for
solving the initial/boundary value problems of the linear
and semilinear heat equations that are fundamental equations
for vast application areas including the diffusion process,
crystal growth, dislocation dynamics, etc. See \cite{brattkus1992siap,brown1989ajm,dargushbanerjee,greengard1990cpam,ibanez2002,jiang2015acom,jiang2017siammms,messner2015sisc,strain1992jcp,strain1994sisc,wang2019jsc}.
For fast Gauss transforms (FGTs) in higher dimensions, see
\cite{beylkin2005jcp,greengard1991fgt,greengard1998nfgt,lee2006dtfgt,sampath2010pfgt,spivak2010sisc,strain1991vfgt,tausch2009sisc,wang2018sisc}.
These FGTs can certainly handle the one-dimensional case, but they
rely on the Hermite and plane wave expansions, complex data structures
and algorithmic steps, leading to nontrivial implementation task.

Here we would like to propose a simple fast algorithm based on
the sum-of-exponentials (SOE) approximations
of the Gaussian kernel:
\be
G(x;\delta)\approx S_n(x;\delta)\coloneqq\sum_{k=1}^n w_k e^{-t_k \frac{|x|}{\sqrt{\delta}}},
\label{soeappr}
\ee
where $w_k$ and $t_k$ are complex weights and nodes of the SOE approximation.
Very often $n$ is chosen to be an even integer and both nodes and weights
appear as conjugate pairs.
Thus, only half of them are needed for computing \cref{gt}.
Since the exponential function
is the eigenfunction of the translation operator, the convolution with the
exponential function as the kernel can be computed via a simple recurrence
in linear time. This leads to an $\mathcal{O}(n(M+N))$ algorithm for computing \cref{gt}.
We would like to remark that SOE approximations have been used extensively
in the design of fast algorithms. See \cite{cheng1999jcp,gimbutas2019fast,greengard1990cpam,jiang2015acom,wang2019jsc}.

We will show that the SOE approximation of the Gaussian kernel converges
geometrically. That is, the maximum error of the approximation
is $\mathcal{O}\left(c^{-n}\right)$.
The optimal constant $c$ in the geometric convergence is difficult to determine.
For up to about ten-digit accuracy, we use
SOE approximations based on best rational approximations
to $e^x$ on the negative real axis $\mathbb{R}^-= (-\infty,0]$
(see, for example, \cite{varga1983,gonchar1987,trefethen2006bit}
and references therein). 
For about ten-digit accuracy, it only requires twelve
exponentials; by symmetry, only six of them are needed in actual
computation. The chief advantage of the current algorithm is that
while achieving comparable performance with existing FGTs in the regime
of practical interests due to the near optimality of the SOE approximations,
its implementation is effortless, requiring only twenty-four lines of
{\sc MATLAB} code. Moreover, 
the algorithm is straightforward for parallelization.

The paper is organized as follows. In \Cref{sec:soeappr}, we discuss
several methods for obtaining SOE approximations of the Gaussian kernel.
The numerical experiments of these methods are presented in \Cref{sec:experiments}.
A simple fast Gauss transform in one dimension based on SOE approximations
is described in \Cref{sec:alg}. Numerical results on the 1D FGT
are presented in \Cref{sec:results} with further discussions
in \Cref{sec:conclusion}.
\section{SOE approximations of the Gaussian kernel}\label{sec:soeappr}
\subsection{Integral representation of the Gaussian}
In \cite{jiang2015acom}, it is shown that the Laplace transform
of the one-dimensional heat kernel $\frac{1}{\sqrt{4\pi t}}e^{-\frac{|x|^2}{4t}}$
is $\frac{1}{2\sqrt{s}}e^{-\sqrt{s}|x|}$. That is
\be
\frac{1}{\sqrt{4\pi t}}e^{-\frac{|x|^2}{4t}}
=\frac{1}{2\pi i}\int_\Gamma e^{st}\frac{1}{2\sqrt{s}}e^{-\sqrt{s}|x|}ds.
\label{heatkernelrep}
\ee
Multiplying both sides by $\sqrt{4\pi t}$, replacing $t$ by $\delta$,
then introducing the change of variables $s=z/\delta$, we obtain
\be
G(x;\delta)=e^{-\frac{x^2}{4\delta}}=
\frac{1}{2\pi i}\int_\Gamma e^{z}\sqrt{\frac{\pi}{z}}
e^{-\frac{\sqrt{z}|x|}{\sqrt{\delta}}}dz.
\label{gaussianrep}
\ee
Here the branch cut of the square root function is chosen to be the negative real axis
$\mathbb{R}^-$ and the branch is chosen to be the principal branch
with $\arg(z)\in (-\pi, \pi]$. We observe that on $\mathbb{C}\setminus \mathbb{R}^-$
the square root function $\sqrt{z}$ always
has positive real part and thus
\be
\left|e^{-\frac{\sqrt{z}|x|}{\sqrt{\delta}}}\right|\le 1
\ee
for all $x\in\mathbb{R}$ and $\delta>0$.
By Cauchy's theorem, the contour $\Gamma$ in \cref{gaussianrep} can be
any contour in the complex plane that starts from $-\infty$ in the
third quadrant, around $0$, and back to $-\infty$ in the second quadrant.
It is clear that the discretization of the integral in \cref{gaussianrep} leads to
an SOE approximation \cref{soeappr} of the Gaussian kernel.
In the following, we will consider integrals of the form
\be
I=\frac{1}{2\pi i}\int_\Gamma e^{z}f(z)dz,
\label{contourintegral}
\ee
and note that for the Gaussian kernel,
\be
f(z)=\sqrt{\frac{\pi}{z}}
e^{-\frac{\sqrt{z}|x|}{\sqrt{\delta}}}.
\label{gf}
\ee

\subsection{Numerical evaluation of the contour integral \cref{contourintegral}}
The numerical evaluation of the contour integral \cref{contourintegral}
has received much attention recently. First, the contour is parametrized 
via $z=z(\theta)$ with $\theta$ a real parameter, leading to
\be
I=\frac{1}{2\pi i}\int_{-\infty}^\infty e^{z(\theta)}f(z(\theta))z'(\theta)d\theta.
\label{contour2}
\ee
The integral in \cref{contour2} is then truncated and discretized 
via either the trapezoidal rule
or the midpoint rule, resulting in an approximation
\be
I_n = \frac{h}{2\pi i }\sum_{k=1}^n e^{z(\theta_k)}f(z(\theta_k))z'(\theta_k),
\label{traprule}
\ee
where $\theta_k$ are $n$ equi-spaced points with spacing $h$.
Substituting \cref{gf} into \cref{traprule}, we obtain an SOE
approximation for the Gaussian kernel
\be
S_n(x;\delta)=\sum_{k=1}^n w_ke^{-t_k\frac{|x|}{\sqrt{\delta}}}
\label{soecontour}
\ee
with
\be
w_k=\frac{h}{2\sqrt{\pi} i }
\frac{z'(\theta_k)}{\sqrt{z(\theta_k)}}e^{z(\theta_k)},
\qquad t_k=\sqrt{z(\theta_k)}.
\ee
Obviously, the efficiency and accuracy of the approximation $I_n$ depends on
the contour and its parametrization. The classical theory summarized in
\cite{stenger1993} 
shows that $I_n$ converges to $I$ subgeometrically with the
rate $\mathcal{O}\left(e^{-c\sqrt{n}}\right)$. Recent developments have improved the
convergence rate to geometric $\mathcal{O}\left(e^{-cn}\right)$ (see, for example,
\cite{trefethen2014sirev} for an excellent review
on applications and analysis of the trapezoidal rule).

There are mainly three kinds of contours. In \cite{talbot1979}, Talbot proposed
cotangent contours; simpler parabolic contours were then proposed in \cite{makarov2000ami};
finally, hyperbolic contours were proposed in \cite{lopez2004anm,lopez2006sinum}.
The hyperbolic contours are applicable to a broader class of functions $f$ in
\cref{contourintegral} in the sense that $f$ is allowed to have singularities
in a sectorial region around $\mathbb{R}^-$.
All three contours have certain parameters that need to be optimized
in order to achieve optimal convergence rate. Much research has been done on the
analysis for the cotangent contours (also
called Talbot contours in literature) (see, for example, \cite{lin2004thesis}).
The optimal choices for all
three contours were first listed in \cite{trefethen2006bit}, with
detailed analysis on the optimization for parabolic and hyperbolic contours appeared in
\cite{weideman2010,weideman2007mcomp}, and the analysis on the modified Talbot contour
in \cite{dingfelder2015na}. We remark that in \cite{weideman2010,weideman2007mcomp}
the total number of quadrature points is $2n$ instead of $n$.

\subsection{Best rational approximations to $e^z$ on $\mathbb{R}^-$}
The parabolic, hyperbolic, and cotangent contours are optimized
within each class. When considering the discrete approximation $I_n$
to $I$, we may generalize $I_n$ and rewrite it as follows:
\be
I_n=-\sum_{k=1}^n c_k f(z_k).
\label{generalapprox}
\ee
And the corresponding SOE approximation for the Gaussian kernel
becomes
\be
S_n(x;\delta)=\sum_{k=1}^n w_ke^{-t_k\frac{|x|}{\sqrt{\delta}}},
\quad w_k=-c_k\sqrt{\frac{\pi}{z_k}}, \quad t_k=\sqrt{z_k}.
\label{bestsoe}
\ee
With a general discrete approximation in the form \cref{generalapprox},
it is natural to ask whether there exists a better approximation
if the points $z_k$ can freely move in a certain region of the complex
plane instead of lying on a curve and the weights $c_k$ are determined
to minimize the quadrature error $I-I_n$. Obviously, this is a very
difficult nonlinear global optimization problem. An interesting connection,
however, is established in \cite{trefethen2006bit} between the approximation
of $I$ by $I_n$ and the best rational approximation of the exponential
function $e^x$ on $\mathbb{R}^-$. More precisely, it is shown there using the
residue theorem and Cauchy's theorem that
\be
I-I_n=\frac{1}{2\pi i}\int_{\Gamma'} \left(e^z-r_{n-1,n}(z)\right)f(z)dz.
\label{errorformula}
\ee
Here $\Gamma'$ is a contour
lying between $\mathbb{R}^-$ and the points $z_k$,
and
\be
r_{n-1,n}(z)=\sum_{k=1}^n \frac{c_k}{z-z_k}
\label{soprep}
\ee
with the index $(n-1,n)$ being
the type of the rational function $r$.
The error formula \cref{errorformula}
indicates that $I-I_n$ would be small if $r_{n-1,n}(z)$ is a good approximation
to $e^z$ on $\mathbb{R}^-$. This is particularly true when all singularities
of $f$ lie on $\mathbb{R}^-$ since $\Gamma'$ can then be chosen arbitrarily
close to $\mathbb{R}^-$, which is the case for the Gaussian kernel.

The latter problem is closely related to the famous ``$1/9$'' problem in
rational approximation theory -- the best rational approximation to
$e^z$ on $\mathbb{R}^-$. The reader is referred to \cite{trefethen2006bit}
for a review on this fascinating problem.
In \cite{gonchar1987}, a rigorous proof showed that
\be
\max_{z\in\mathbb{}R^-}\left|e^z-r^\ast_{n,n}(z)\right| = \mathcal{O}\left(9.28903\cdots^{-n}\right),
\ee
where $r^\ast_{n,n}$ is the best rational approximation of type $(n,n)$;
in \cite{varga1983},
a numerical algorithm based on the classical Remez algorithm was devised to
calculate the best rational approximation $r^\ast_{n,n}$ and the approximation
error to very high precision;
in \cite{trefethen2006bit}, 
an algorithm based on the Carath{\'e}odory--Fej{\' e}r (CF) method was developed
to compute a nearly best rational approximation of type $(n-1,n)$
for $n$ up to $14$ (see the {\sc MATLAB}
code cf.m in \cite{trefethen2006bit}).

\subsection{Global reduction via balanced truncation method for SOE approximations}
The type $(n-1,n)$ rational function $r_{n-1,n}$ in \cref{soprep}
is also called sum-of-poles (SOP) function. When all poles lie in, say, the left half
of the complex plane, the balanced truncation method can be used to
obtain an SOP approximation of smaller number of poles with a prescribed
$L^{\infty}$ error. In \cite{greengard2018sisc}, it was observed that the Laplace
transform of an exponential function is a pole function, and thus the balanced truncation
method can be used to reduce the number of exponentials in an SOE approximation as well.
There are extensive theoretical investigations on infinite Hankel matrices
\cite{aak1968a,aak1968b,aak1971,peller2003}, and related algorithmic developments
on the balanced truncation method \cite{glover1984ijc}, the CF method for rational
approximations \cite{gutknecht1983b,gutknecht1984jat,gutknecht1982sinum,gutknecht1983a,trefethen1983sinum,trefethen1983tams}, and Prony's method for SOE approximations
\cite{beylkin2005,beylkin2010}. Here we will use a simplified balanced truncation method
outlined in \cite{xu2013jsc} to try to further reduce the number of exponentials
for a given SOE approximation of the Gaussian kernel.
\section{Numerical Experiments on SOE approximations of the Gaussian kernel}
\label{sec:experiments}
We now present numerical experiments on SOE approximations of the Gaussian
kernel using methods discussed in \Cref{sec:soeappr}. The maximum error
\be
E_n=\max_{x\in \mathbb{R}}|G(x;\delta)-S_n(x;\delta)|
\label{maxerror}
\ee
is estimated by setting
$\delta=1$ and sampling $x$ at $0$ and $100,000$ logarithmically equally spaced
points on $[10^{-5},10^2]$ (note that $E_n$ is independent of $\delta$ for any $\delta>0$).
\subsection{Optimal parabolic, hyperbolic, and modified Talbot contours}
\begin{figure}[!ht]
\centering
\includegraphics[height=32mm]{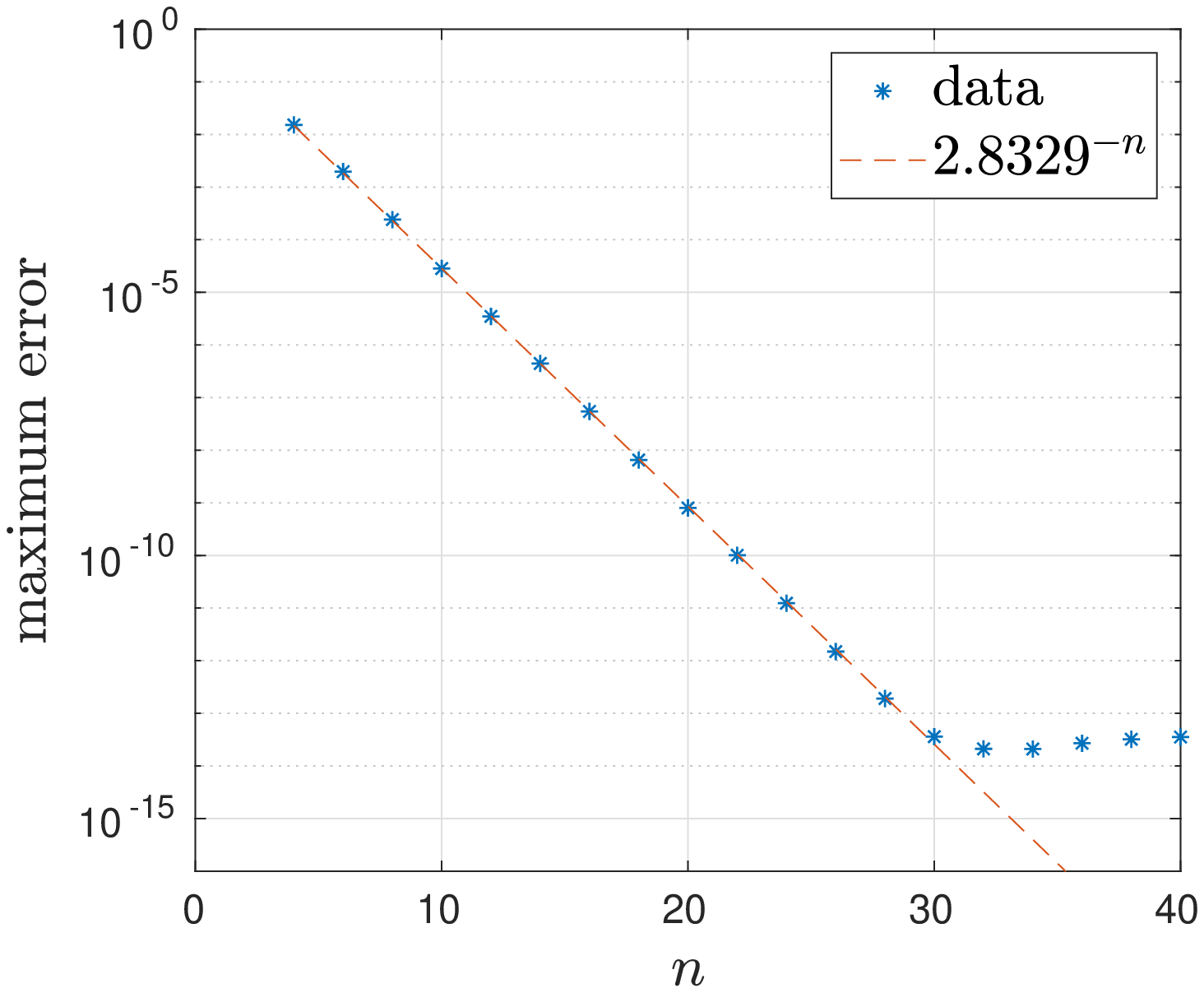}
\includegraphics[height=32mm]{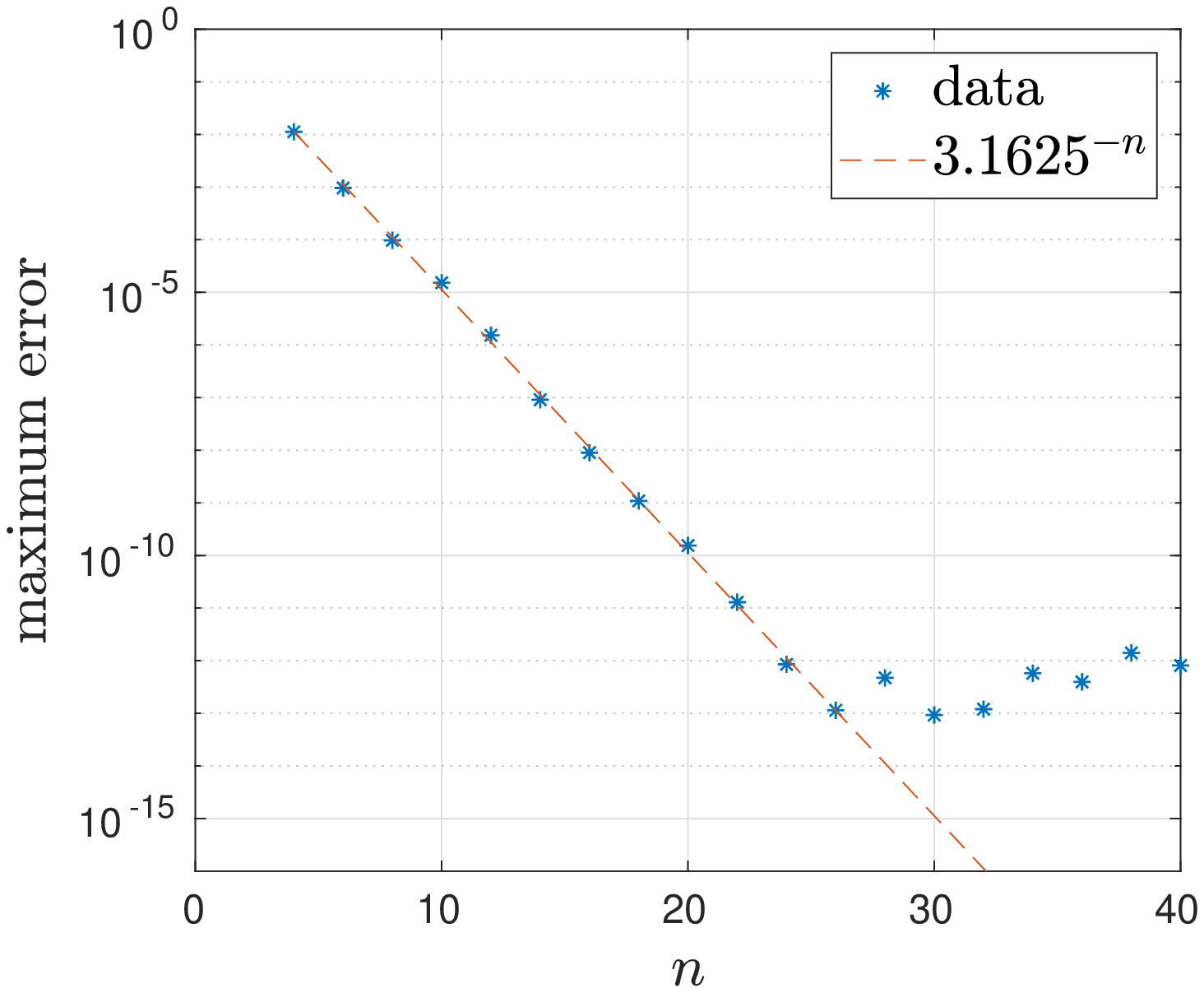}
\includegraphics[height=32mm]{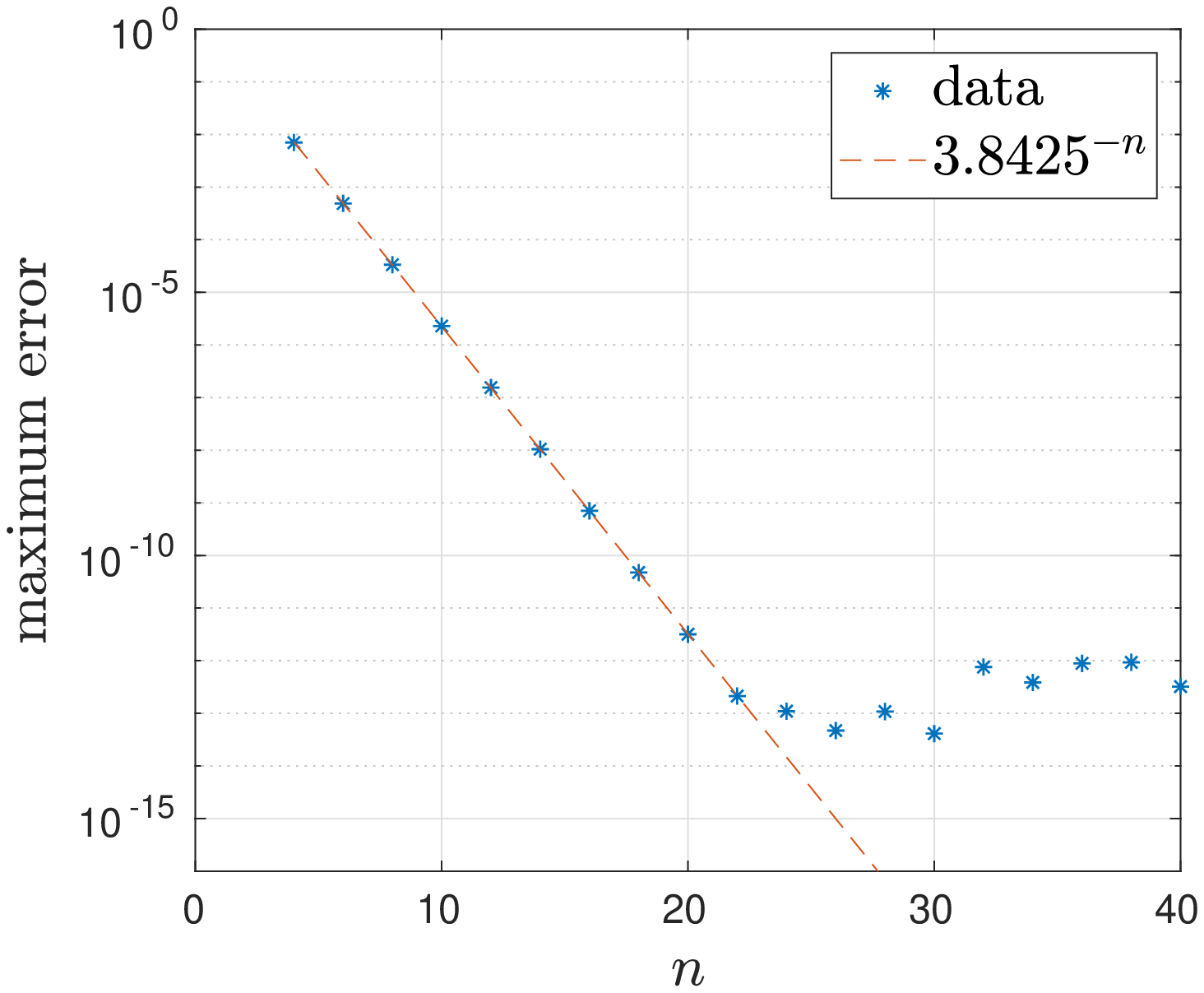}
\caption{Maximum error of the SOE approximation of the Gaussian kernel
  as a function of $n$.
  Left: optimal parabolic contour. Middle: optimal hyperbolic contour.
  Right: optimal modified Talbot contour. See \cite{trefethen2006bit}.
  Dashed lines show the estimated convergence rate via least squares fitting
  of the data points.
}
\label{fig1}
\end{figure}

\Cref{fig1} shows $E_n$ as a function of $n$ when
the contours in \cite{trefethen2006bit} are discretized via the midpoint rule.
We observe that convergence rates are very close to theoretical values
in \cite{trefethen2006bit} for all three contours.
The maximum errors saturate at about 13-digit
accuracy due to roundoff errors.

\subsection{Further reduction by the balanced truncation method}
\label{sec:reduction}
\begin{figure}[!ht]
\centering
\includegraphics[height=40mm]{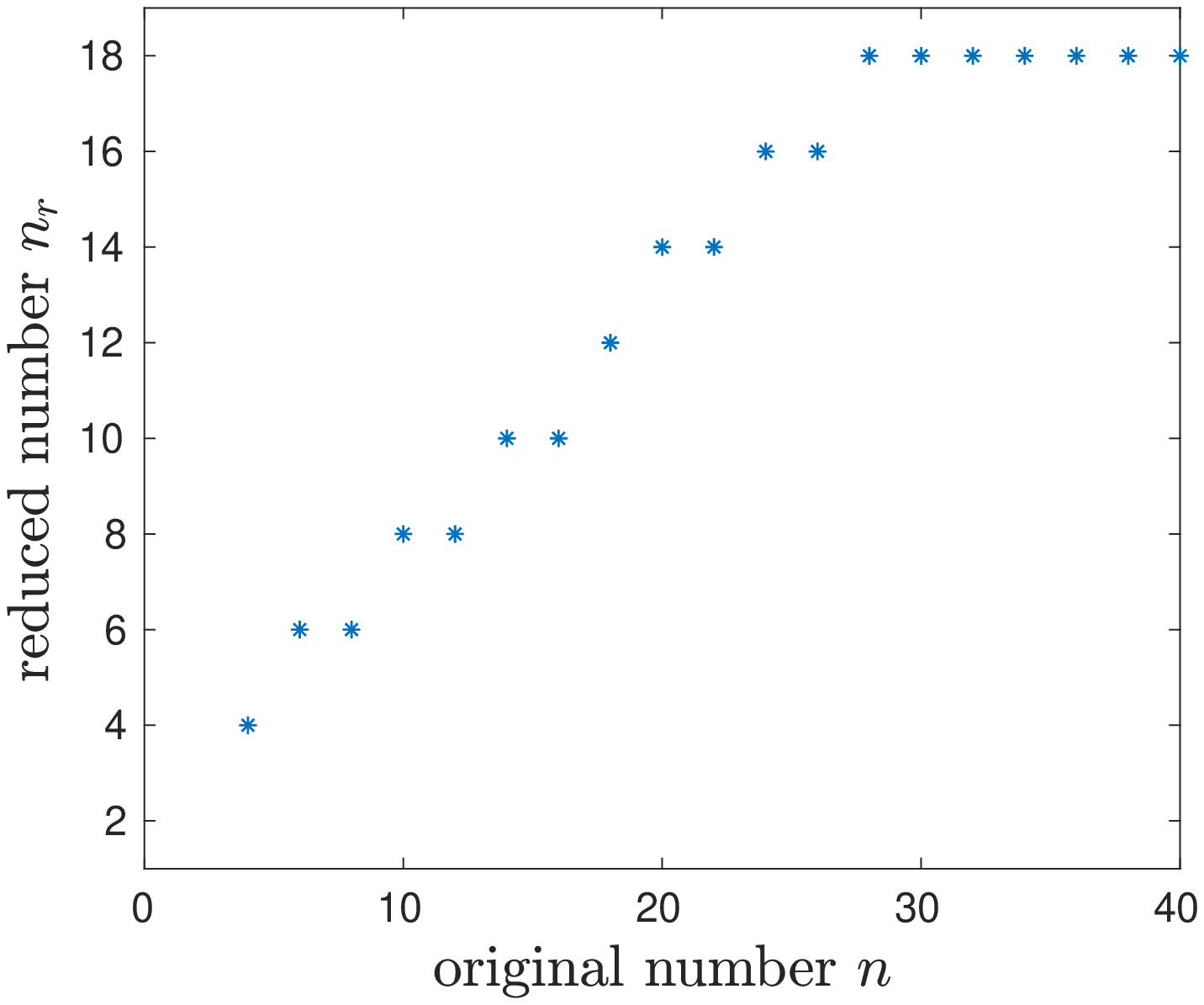}
\includegraphics[height=40mm]{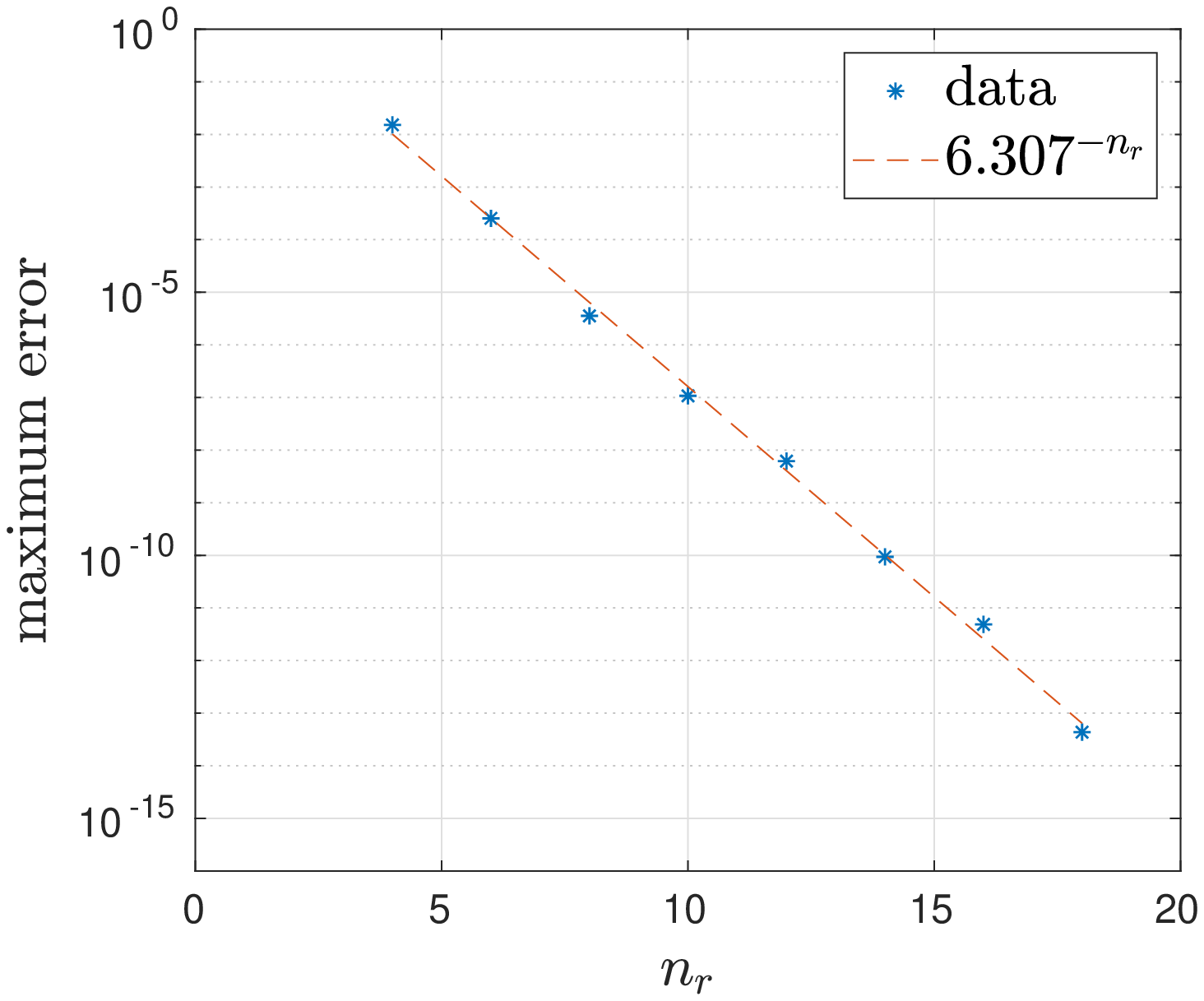}

\vspace{4mm}

\includegraphics[height=40mm]{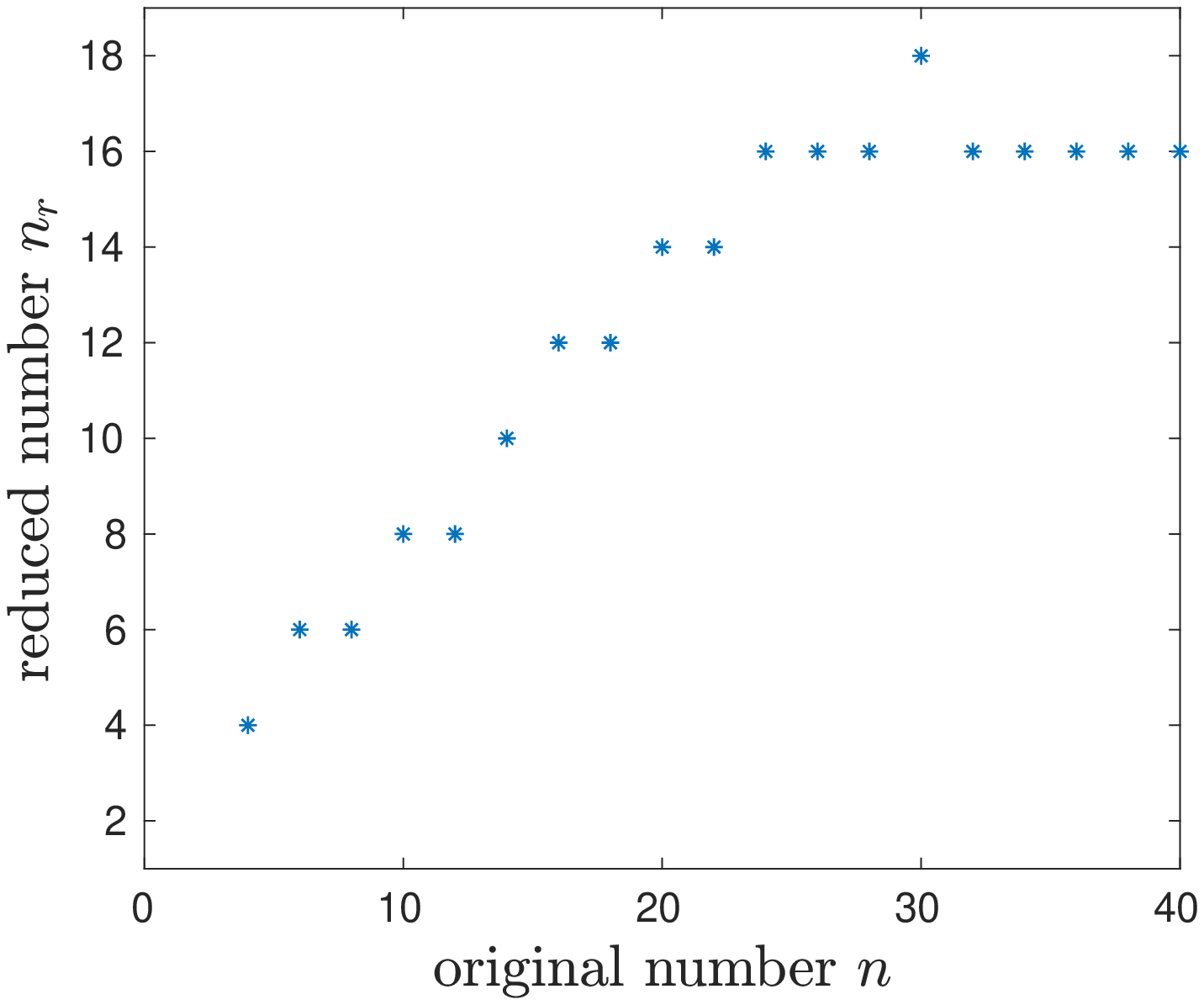}
\includegraphics[height=40mm]{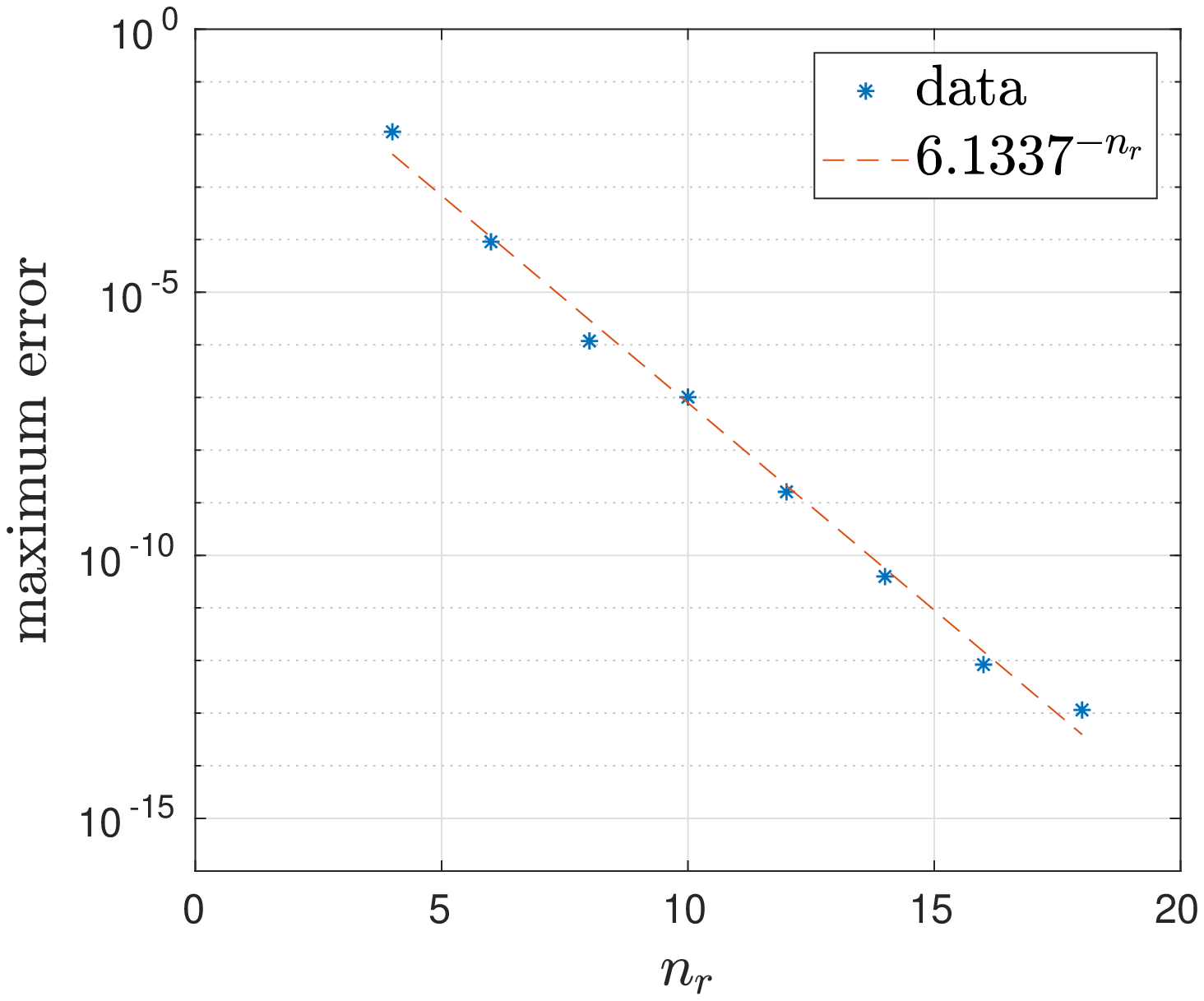}

\vspace{4mm}

\includegraphics[height=40mm]{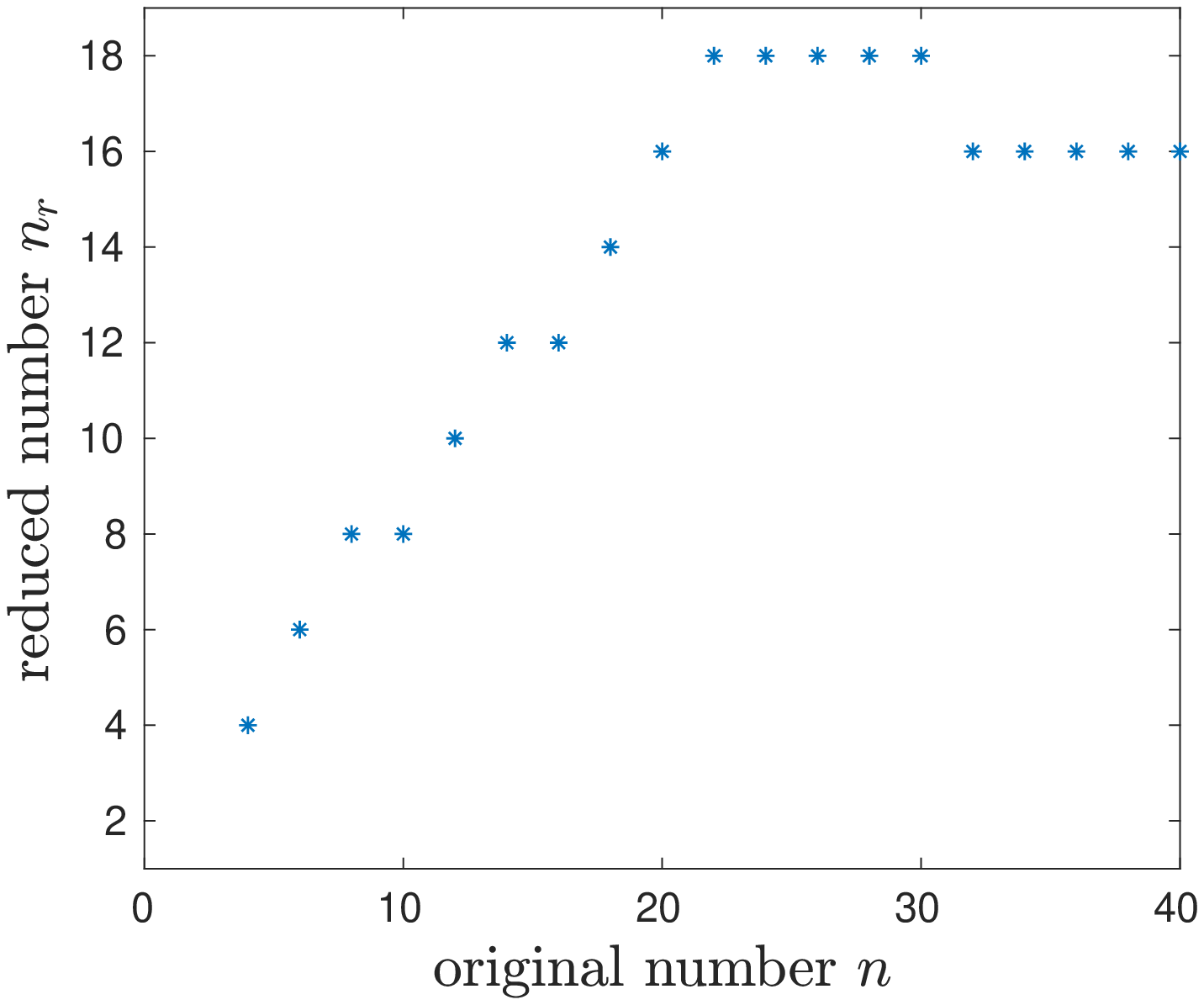}
\includegraphics[height=40mm]{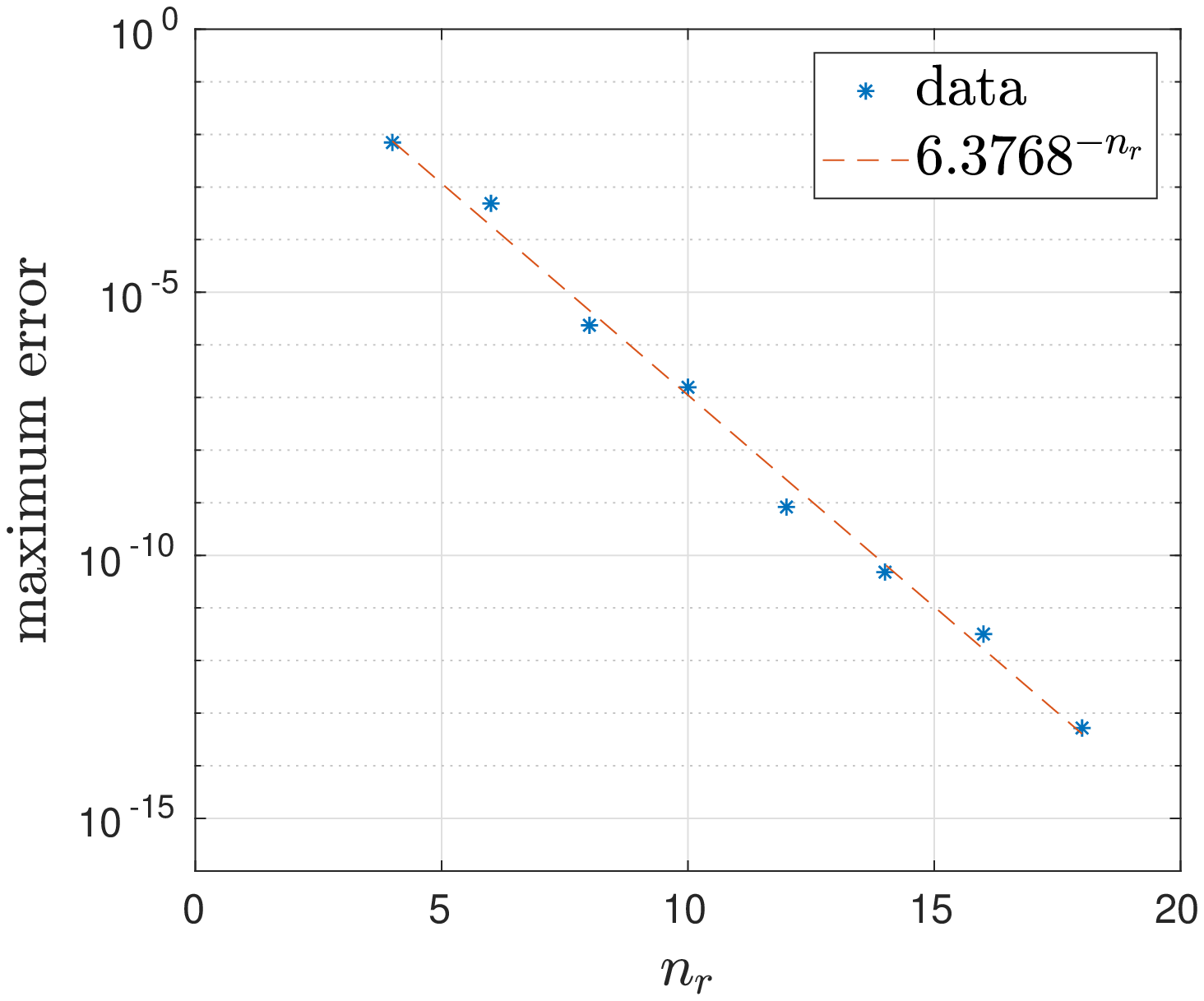}

\caption{Effect of the balanced truncation method to reduce the number of exponentials
  in the SOE approximation of the Gaussian kernel. Left: reduced number $n_r$ obtained
  by the balanced truncated method 
  versus original number $n$ obtained by the midpoint rule discretization of the
  contour integrals. Right: estimated convergence rate of the reduced SOE approximation.
  Top: optimal parabolic contour. Middle: optimal hyperbolic contour.
  Bottom: optimal modified Talbot contour.}
\label{fig2}
\end{figure}

We then try to use the balanced truncation method to reduce the number of exponentials
in SOE approximations. The prescribed precision for the balanced truncation
method is set to $E_n/3$ so that the reduced SOE approximation has about the same
accuracy as the original one. It turns out that SOE approximations from these
three optimal contours can all be reduced. The numerical results are summarized in
\Cref{fig2}. We observe that for all three contours, the reduced number $n_r$ of
exponentials in the SOE approximation saturates at $18$, and the convergence rate
is improved to about $\mathcal{O}\left(6.3^{-n}\right)$.

\subsection{Best rational approximations}
\begin{figure}[!ht]
\centering
\includegraphics[height=40mm]{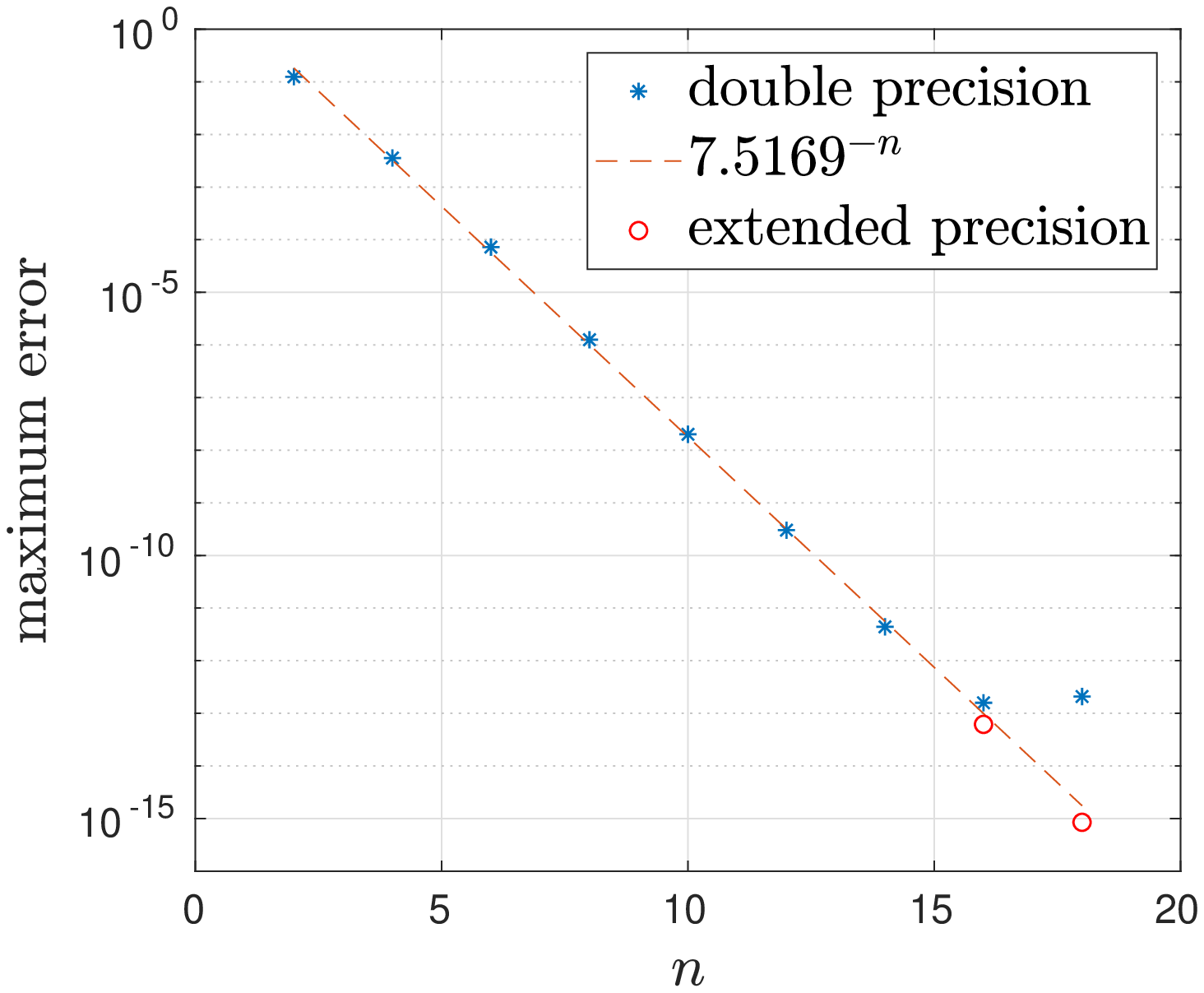}
\includegraphics[height=40mm]{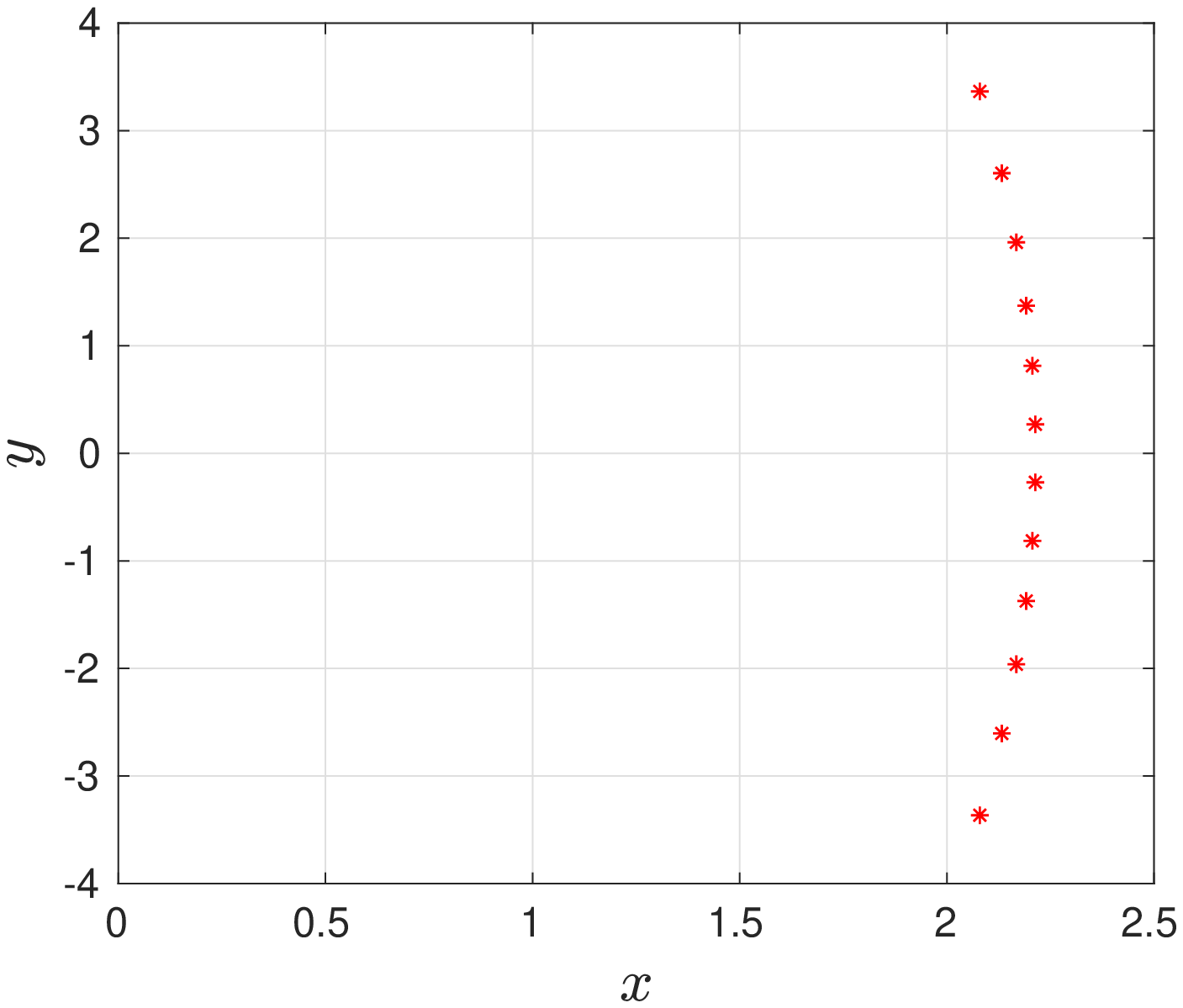}
\caption{The SOE approximation of the Gaussian kernel obtained by
  the best rational approximation to the exponential function on the
  negative real axis.
  Left: maximum error $E_n$ as a function of $n$. Dashed line is the least
  squares fit to the first seven data points. The red circles for $n=16, 18$
  are calculated in extended precision. Right: locations of exponential nodes
  $t_k$ for $n=12$.
}
\label{fig3}
\end{figure}

We now turn to SOE approximations obtained by best rational
approximations to the exponential function on the negative real axis.
For $n\le 14$, we use the {\sc MATLAB} code cf.m in \cite{trefethen2006bit}
to compute pole locations and weights. For $n=16$ and $18$, we have
implemented the CF algorithm in cf.m in {\sc Fortran} and run
the code in quadruple precision to obtain pole locations and weights. The attempt
of using the balanced truncation method on these SOE approximations does not
lead to any further reduction.
\Cref{fig3} summarizes numerical results. We observe that the convergence rate
is about $\mathcal{O}\left(7.5^{-n}\right)$, which is better than those obtained
by optimal contours even after further reduction step by the balanced truncation
method. Second, the roundoff error makes the approximation saturate
at about 13-digit accuracy in double precision arithmetic. By symmetry, only half
of the exponentials are needed in actual computation when $n$ is even. Hence, we
only need $3$ exponentials to achieve four-digit accuracy, and $6$ exponentials
to achieve about ten-digit accuracy.

\subsection{Effect of roundoff errors}
\label{sec:roundoff}
We observe that all SOE approximations discussed in the previous subsections
have errors saturated at about 13-digit. This is due to the ill-conditioning of
the inverse Laplace transform and similar phenomena have been observed
many times (see, for example, \cite{lopez2006sinum,weideman2007mcomp}).
Efforts have been made in stablizing the inverse Laplace transforms. In
\cite{lopez2006sinum}, a detailed analysis on the roundoff error for hyperbolic
contours is presented and an additional parameter $\theta$ is introduced to balance
the roundoff error with other errors (i.e., the discretization error and the truncation
error). In \cite{weideman2010}, a strategy is proposed to stabilize the parabolic
contours. Finally, in \cite{dingfelder2015na}, the stabilization of the modified
Talbot contours is studied.

We have tested the stabilizing techniques presented in these papers.
All of them do stabilize the SOE approximations for the Gaussian kernel.
However, the parabolic and Talbot contours have the approximation error
saturated at about 14-digit. The accuracy of the SOE approximation using
unoptimized hyperbolic contours in \cite{lopez2006sinum} can go down to
almost machine precision in double precision arithmetic. And the balanced truncation
method will reduce the number of exponentials and improve the convergence
rate greatly.

\begin{figure}[!ht]
\centering
\includegraphics[height=32mm]{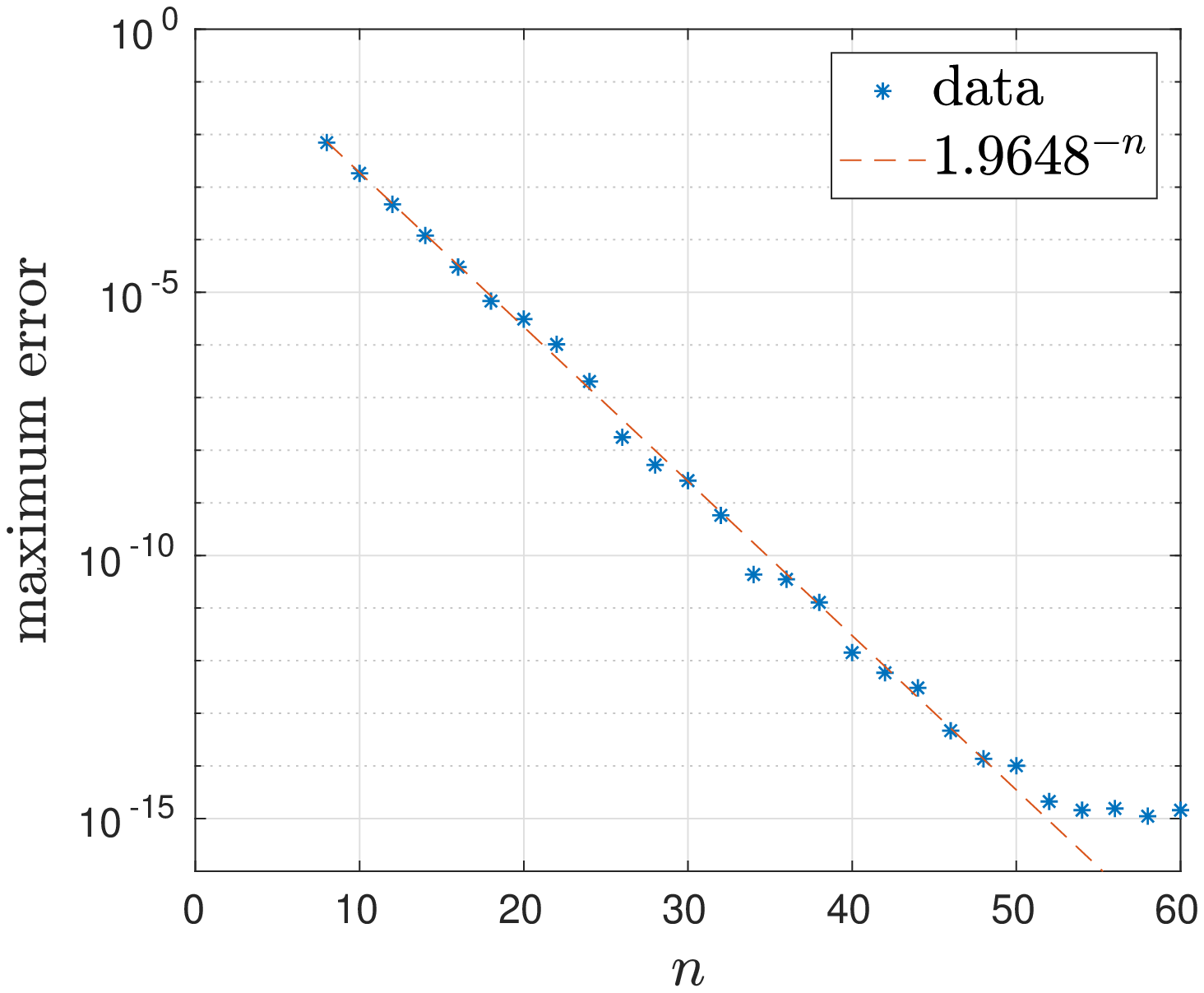}
\includegraphics[height=32mm]{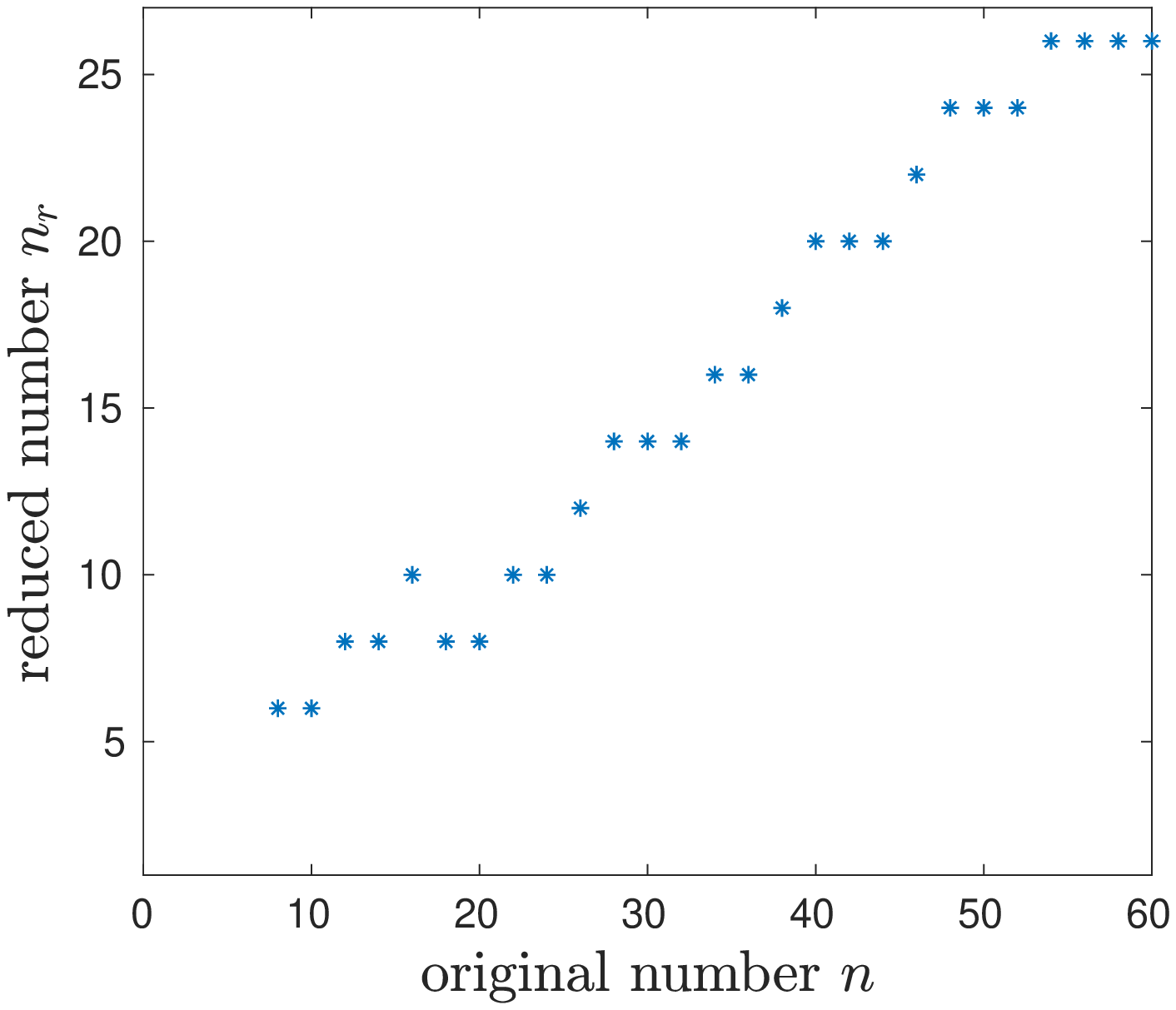}
\includegraphics[height=32mm]{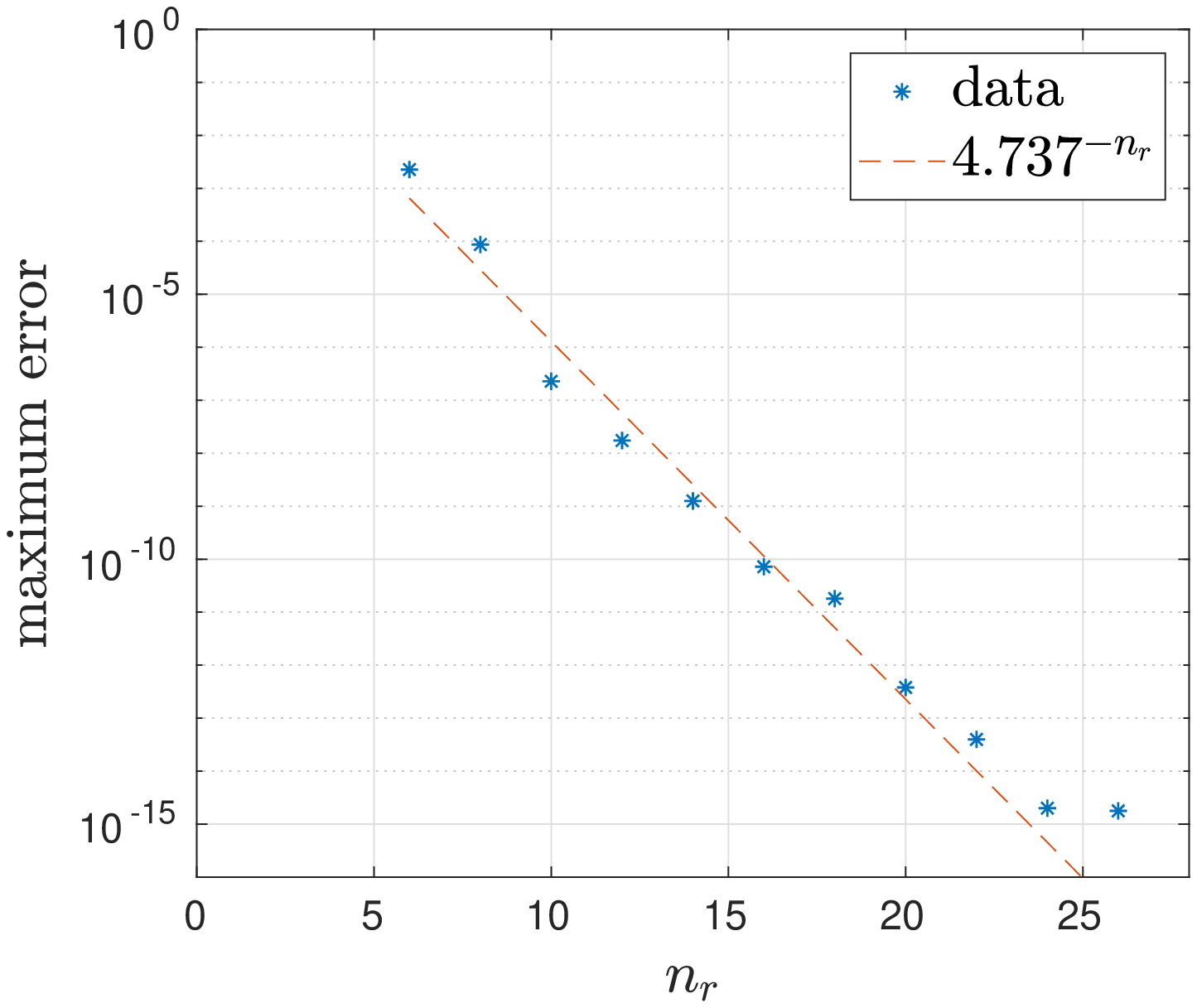}
\caption{SOE approximation of the Gaussian kernel with smaller roundoff error.
  Left: errors of SOE approximations by a
  hyperbolic contour in \cite{lopez2006sinum}. Middle: reduction by
  the balanced truncation method. Right: improved convergence rate. 
  Dashed lines show the estimated convergence rate via least squares fitting
  of the data points.
}
\label{fig4}
\end{figure}

In a certain sense, the ill-conditioning of the inverse Laplace transform 
can be viewed as an advantage when we try to construct SOE approximations for
the Gaussian kernel. To be more precise, there are a wide range of parameters
and contours that will lead to quite different SOE approximations with about same
accuracy and the same number of exponentials. Such phenomena are already shown
in \Cref{sec:reduction}. To demonstrate that the SOE approximation
can achieve close to machine precision, we show one of our numerical experiments
in \cref{fig4}. Here we use the hyperbolic contour in \cite{lopez2006sinum}
$z=\lambda(1-\sin(0.8+iu))$
with $a=\cosh^{-1}(2/((1-\theta)\sin(0.8)))$, $\lambda=0.6\pi(1-\theta)n/a$,
the step size $h=2a/n$. The parameter $\theta$ is set to $1/4$ for $n\le 16$
and $12/n$ otherwise. We observe that at $n_r=24$, the error of the reduced
SOE approximation is about $2.0\times 10^{-15}$, at least ten times smaller
than those by other methods that we studied.

\section{A simple fast Gauss transform in one dimension}\label{sec:alg}
We now present a fast algorithm for the Gauss transform in one dimension
using SOE approximations. The algorithm has been known in the field of
fast algorithms. Most recently it has been used in \cite{gimbutas2019fast},
where SOE approximations are developed for the kernel $\frac{1}{x}$ on $[\delta,R]$.
Since the SOE approximations for the Gaussian kernel is valid in the whole
real axis, the algorithm in \cite{gimbutas2019fast} can be simplified a bit,
which is discussed briefly below.

We first sort both target and source points in ascending order. In the following
discussion, we will assume that targets and sources are already sorted.
Substituting the SOE approximation \cref{soeappr} into \cref{gt} and exchange
the order of summations, we obtain
\be
u_i\approx\sum_{k=1}^nw_k h_{k,i}, \qquad i=1,\ldots, M
\label{soegt}
\ee
where
\be
h_{k,i}=\sum_{j=1}^N\alpha_je^{-t_k\frac{|x_i-y_j|}{\sqrt{\delta}}}.
\label{hmode}
\ee
We further split $h_{ki}$ into two terms
\be
h_{k,i}=h_{k,i}^++h_{k,i}^-,
\label{split}
\ee
where
\be
h_{k,i}^+=\sum_{y_j\le x_i} \alpha_je^{-t_k\frac{x_i-y_j}{\sqrt{\delta}}},\qquad
h_{k,i}^-=\sum_{y_j> x_i} \alpha_je^{-t_k\frac{y_j-x_i}{\sqrt{\delta}}}.
\label{split2}
\ee
It is clear that $h_{k,i}^+$ satifies the forward recurrence
\be
h_{k,i+1}^+=e^{-t_k\frac{x_{i+1}-x_i}{\sqrt{\delta}}}h_{k,i}+
\sum_{x_i<y_j\le x_{i+1}} \alpha_je^{-t_k\frac{x_{i+1}-y_j}{\sqrt{\delta}}},
\label{frecurrence}
\ee
and $h_{k,i}^-$ satifies the backward recurrence
\be
h_{k,i-1}^-=e^{-t_k\frac{x_{i}-x_{i-1}}{\sqrt{\delta}}}h_{k,i}+
\sum_{x_i\ge y_j> x_{i-1}} \alpha_je^{-t_k\frac{y_j-x_{i-1}}{\sqrt{\delta}}}.
\label{brecurrence}
\ee
Hence, $h_{k,i}^\pm$ can be computed in $\mathcal{O}(N+M)$ time
for each $k$ and all $i=1,\ldots,M$.

Instead of presenting the pseudocode for the whole algorithm, we
show 24-line {\sc MATLAB} code for the case when the target points
are identical to source points in \cref{fig5}. Several remarks are as
follows. First, in order to save some space
we have called cf.m in \cite{trefethen2006bit} directly in the code.
Slightly better implementation may precompute and store nodes and weights
of the SOE approximation. Second, due to the absence of a genuine compiler
in {\sc MATLAB}, we have flipped the backward pass into a forward loop
for faster memory access. The code runs at about one million
points per second on a laptop with Intel(R) 2.10GHz i7-4600U CPU. Finally,
when the target points and source points are different, the forward pass
is replaced by the code fragment in \cref{fig6}. The changes to the backward
pass are similar and thus omitted. 

\begin{figure}[!ht]
\centering
\begin{lstlisting}
  function u=fgt1d(x,alpha,delta)
  [zs,cs] = cf(12);                % call cf.m in Trefethen et al. (2006) 
  zs = zs(1:2:12); cs = cs(1:2:12);% take only half of the exponentials
  ws = -2*sqrt(pi)*cs./sqrt(zs);   % SOE approximation weights
  ts = sqrt(zs/delta);             % SOE approximation nodes
  [xs,I] = sort(x);                % sort the points in ascending order
  beta = alpha(I);                 % align the strength vector with points
  e1 = exp(-ts*diff(xs));          % compute all needed complex exponentials
  n2 = length(ts); nx = length(x);
  hp = zeros(n2,nx); 
  e2 = ones(n2,1); 
  hp(:,1) = beta(1)*e2;
  for i = 2:nx                     % forward pass
    hp(:,i) = beta(i)*e2+e1(:,i-1).*hp(:,i-1);
  end
  hm = zeros(n2,nx); 
  e1 = fliplr(e1); 
  beta = fliplr(beta);
  for i = 2:nx                     % backward pass
    hm(:,i) = e1(:,i-1).*(beta(i-1)+hm(:,i-1));
  end
  u = real(ws.'*(hp+fliplr(hm)));  % sum over all exponential modes
  Iinv(I) = 1:nx;
  u = u(Iinv);                     % reorder u by the original target order
\end{lstlisting}
\caption{{\sc MATLAB} code for the fast Gauss transform in one dimension when the
  target points are identical to the source points. 
}
\label{fig5}
\end{figure}
\begin{figure}[!ht]
\centering
\begin{lstlisting}
  nx = length(x); ny = length(y);
  [xs,I] = sort(x);                % sort the target points
  [ys,J] = sort(y);                % sort the source points
  beta = alpha(J);                 % align the strength vector with sources
  e1=exp(-ts*diff(xs));            % compute needed exponentials
  hp = zeros(length(ts),nx);
  i = 1; j = 1;
  while i <= nx && j <= ny         % sweep through targets and sources 
    if ys(j) < xs(i)               % contribution from different sources
      hp(:,i) = hp(:,i)+beta(j)*exp(-ts*(xs(i)-ys(j)));
      j = j+1;
    elseif ys(j) == xs(i)          % contribution from the identical source
      hp(:,i) = hp(:,i)+beta(j);
      j = j+1;
    else                           % forward recurrence over targets
      if i < nx
        hp(:,i+1) = e1(:,i).*hp(:,i);
      end
      i = i+1; 
    end
  end
\end{lstlisting}
\caption{{\sc MATLAB} code fragment for the forward pass when the target points and source
  points are different.
}
\label{fig6}
\end{figure}

\section{Numerical results}\label{sec:results}
We have implemented the algorithms in \Cref{sec:alg} in Fortran. As to the sorting,
we modify the function dlasrt.f in LAPACK 3.8.0 so that it outputs the
sorted array $xs$ and an integer
arrary $I$ with $xs=x(I)$. The function uses Quick Sort, reverting to Insertion sort
on arrays of size $\le 20$. Thus, the average complexity of the sorting step in
our implementation is $\mathcal{O}(N\log N)$. 
The code is compiled using gfortran 6.3.0 with -O3 option. The results
shown in this section were obtained on a laptop with Intel(R) 2.10GHz i7-4600U CPU
in a single core mode.
In some applications, one may need to apply the Gauss transform
many times with fixed locations of points but different strength vector $\alpha$.
In this case, it is advantageous to pre-sort the points and precompute
all needed complex exponentials and store them in a table.

We first test the case when the target points are identical to the source
points. We have tested the uniform distribution (on $[0,1]$) case and Chebyshev points case.
\Cref{table1} shows the results for one sample run of uniform distribution case
with $\delta=1$. In the table,
$N$ is the total number of points in the Gauss transform. $n_{\rm e}=n/2$ is the actual
number of complex exponentials needed in the computation. $t_{\rm sort}$ is the time
for the sorting step. $t_{\rm pre}$ is the time for precomputing complex exponentials.
$t_{\rm rem}$ is the time for the remaining calculations. $t_{\rm total}$ is the total
computational time. The time is measured in seconds.
And the error is the estimated maximum relative error computed
at $100$ randomly selected target points. The results for Chebyshev points are
similar except that the sorting step takes much less time and thus omitted.
From the table, we observe that the most expensive step is computing all complex
exponentials, and the sorting step already takes a significant portion in the
total computational time for the numbers $N$ shown here. Indeed, if we pre-sort
the points and precompute the exponentials, the code would run as fast as Quick Sort.
\begin{table}[!ht]
  \caption{Results of the FGT on identical target and source points with uniform
    distribution.
  }
  \sisetup{
%  scientific-notation = fixed,    
%  fixed-exponent = 0,
  tight-spacing=true
  }
%  table-auto-round=true,
%  table-format = 1.1e-2,
\centering
\begin{tabular}{S[group-separator = {,},table-format=8.0e0]
    S[table-format=4.0e0]
    S[scientific-notation = fixed,fixed-exponent = 0]
    S[scientific-notation = fixed,fixed-exponent = 0]
    S[scientific-notation = fixed,fixed-exponent = 0]
    S[scientific-notation = fixed,fixed-exponent = 0]
    S[scientific-notation = true,table-format=3.1e-2]}
  \toprule
  ${N}$  & ${\,\,\,\,\,\,\,\,\,\,n_{\rm e}}$ & ${t_{\rm sort}}$ & ${t_{\rm pre}}$ & ${t_{\rm rem}}$ & ${t_{\rm total}}$ & {Error}  \\
  \midrule
    100000 & 3 & 0.8D-02 & 0.16D-01 & 0.8D-02 & 0.36D-01 & 0.44D-05 \\
   1000000 & 3 & 0.92D-01 & 0.16D+00 & 0.68D-01 & 0.34D+00 & 0.43D-05 \\
   10000000 & 3 & 0.10D+01 & 0.16D+01 & 0.64D+00 & 0.38D+01 & 0.43D-05 \\
   &&&&&&\\
    100000 & 4 & 0.8D-02 & 0.20D-01 & 0.4D-02 & 0.32D-01 & 0.55D-07 \\
   1000000 & 4 & 0.88D-01 & 0.21D+00 & 0.80D-01 & 0.40D+00 & 0.55D-07 \\
   10000000 & 4 & 0.10D+01 & 0.21D+01 & 0.78D+00 & 0.44D+01 & 0.55D-07 \\
   &&&&&&\\
    100000 & 5 & 0.8D-02 & 0.24D-01 & 0.4D-02 & 0.44D-01 & 0.63D-09 \\
   1000000 & 5 & 0.92D-01 & 0.26D+00 & 0.96D-01 & 0.46D+00 & 0.62D-09 \\
   10000000 & 5 & 0.10D+01 & 0.25D+01 & 0.94D+00 & 0.51D+01 & 0.56D-09 \\
   &&&&&&\\
    100000 & 6 & 0.8D-02 & 0.32D-01 & 0.8D-02 & 0.48D-01 & 0.76D-11 \\
   1000000 & 6 & 0.88D-01 & 0.30D+00 & 0.12D+00 & 0.53D+00 & 0.49D-11 \\
  10000000 & 6 & 0.10D+01 & 0.31D+01 & 0.11D+01 & 0.57D+01 & 0.95D-10 \\
 \bottomrule
 \end{tabular}
\label{table1}
\end{table}

\begin{figure}[!ht]
\centering
\includegraphics[height=40mm]{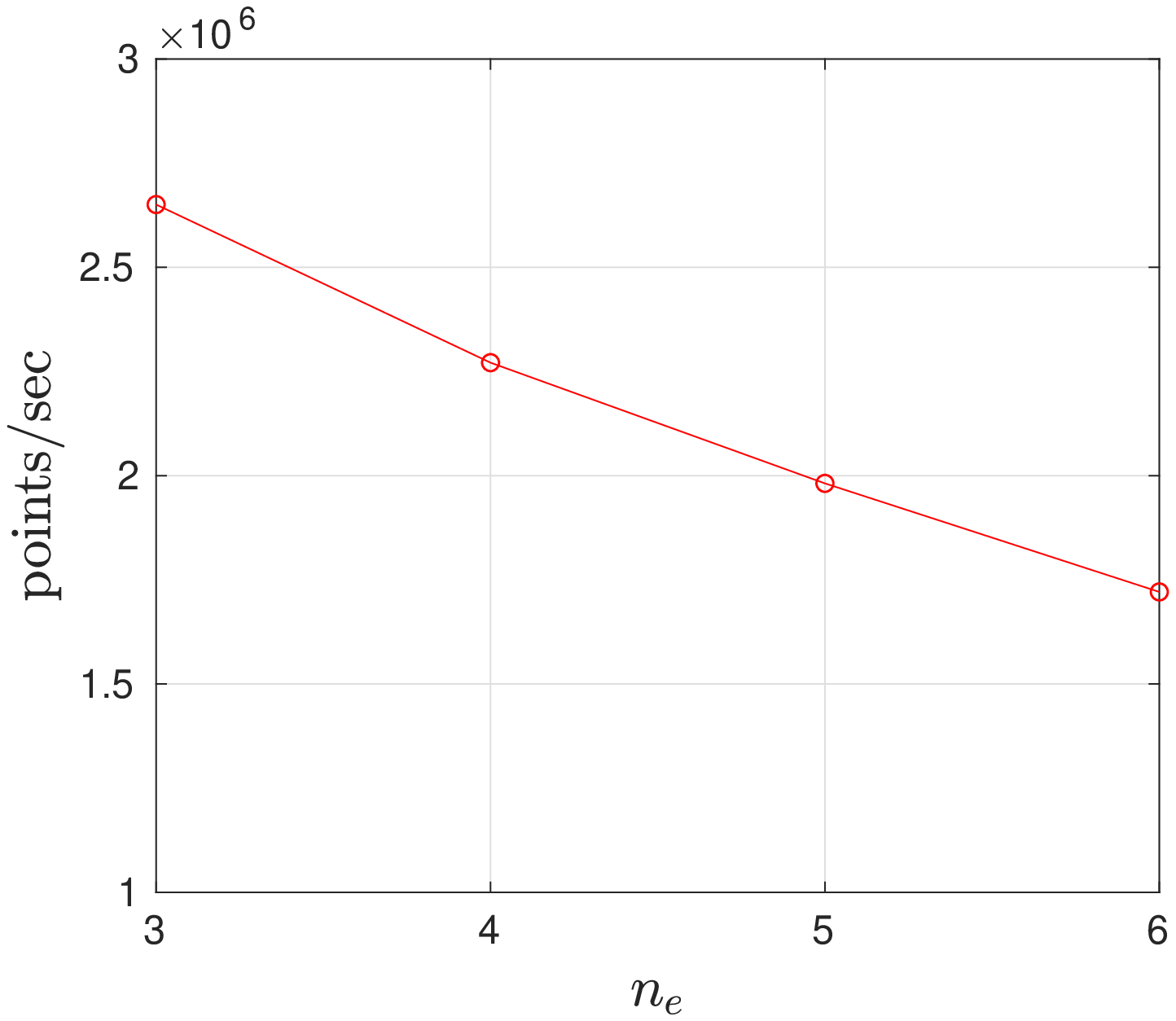}
\includegraphics[height=40mm]{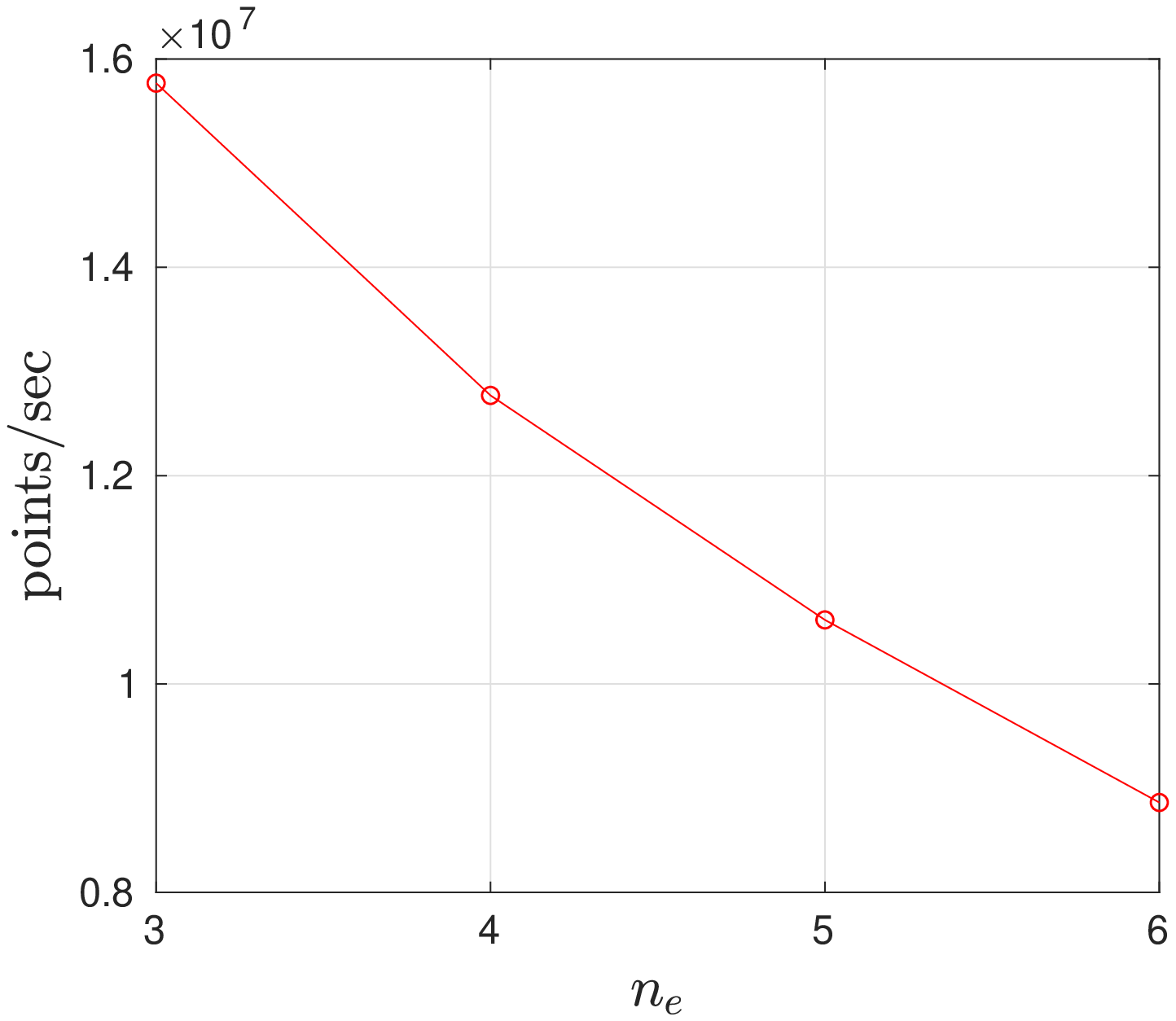}
\caption{Throughput of the 1D FGT as a function of $n_e$
  when targets are identical to sources. $n_e$ is half of the number of exponentials
  in the SOE approximation of the Gaussian kernel.
  Left: Throughput of the FGT without any precomputation.
  Right: Throughput of the FGT when points are pre-sorted and exponentials
  are precomputed and stored.
}
\label{fig7}
\end{figure}

The throughput is shown in \cref{fig7}, where each data point is obtained
by averaging $10$ sample runs with $N$ ranging from $10^6$ to $10^7$.
We observe that the throughput of the whole algorithm ranges from $2.6$ to $1.7$ million
points per second when $n_e$ increases from $3$ to $6$. If we pre-sort
the points and precompute all complex exponentials, the throughput increases
to $16$ to $9$ million points per second.

\begin{figure}[!ht]
\centering
\includegraphics[height=45mm]{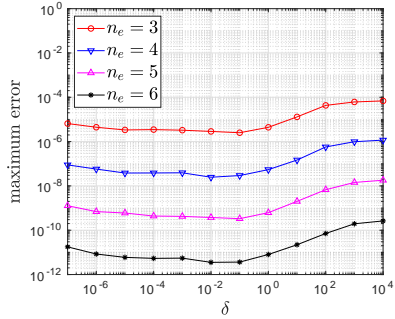}
\caption{Accuracy of the FGT as a function of $\delta$.
}
\label{fig8}
\end{figure}

It is clear that the complexity of the algorithm is independent of $\delta$.
\Cref{fig8} shows the accuracy dependence on $\delta$ for $\delta=10^{-7},\ldots,10^4$.
We observe that the errors slowly increases as $\delta$ increases, and
saturates at the errors shown in \cref{fig3}. This is because SOE approximations
obtained by best rational approximations have the largest error at or near the origin,
and all Gaussians will cluster around the origin as $\delta$ increases.

\begin{figure}[!ht]
\centering
\includegraphics[height=42mm]{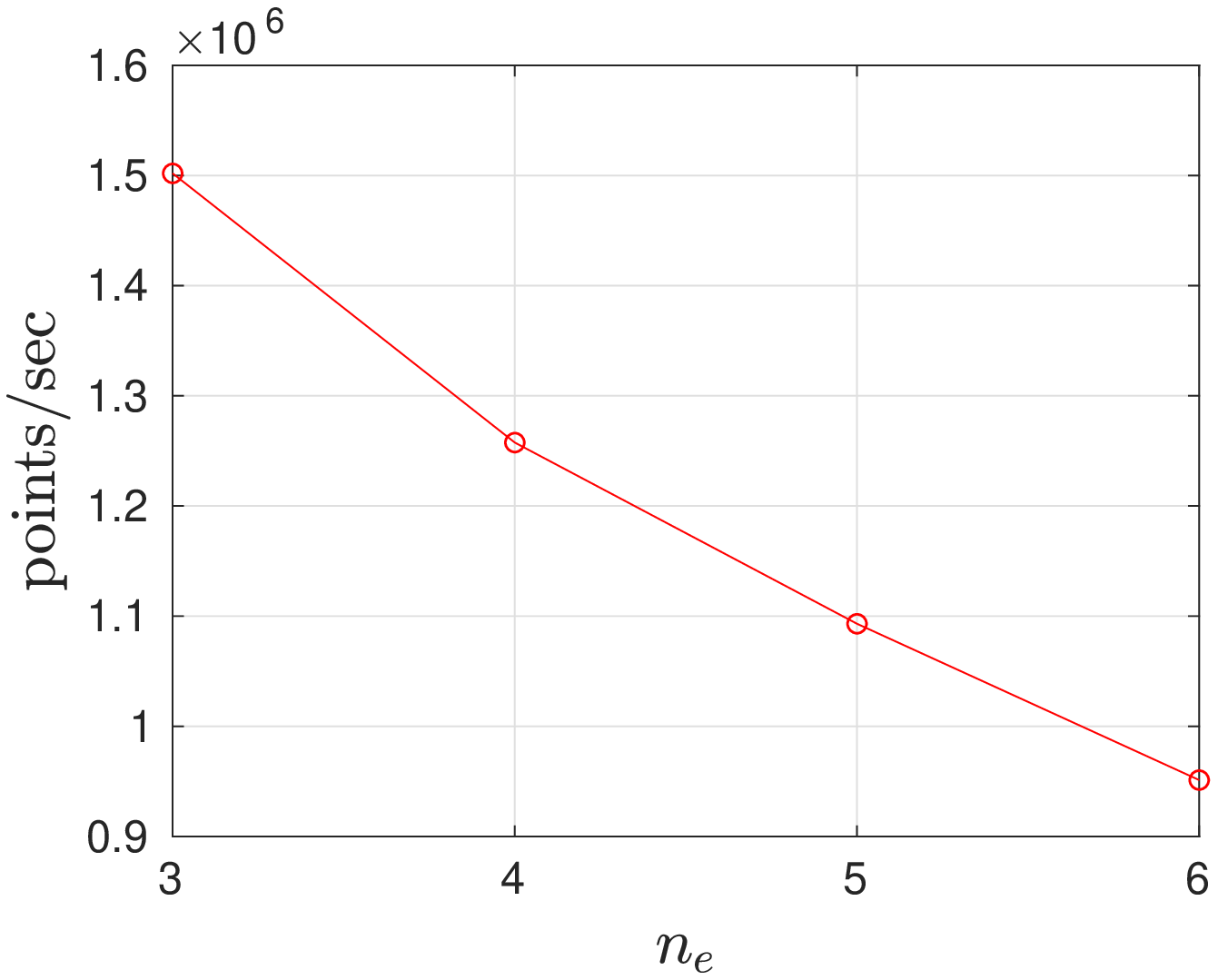}
\includegraphics[height=42mm]{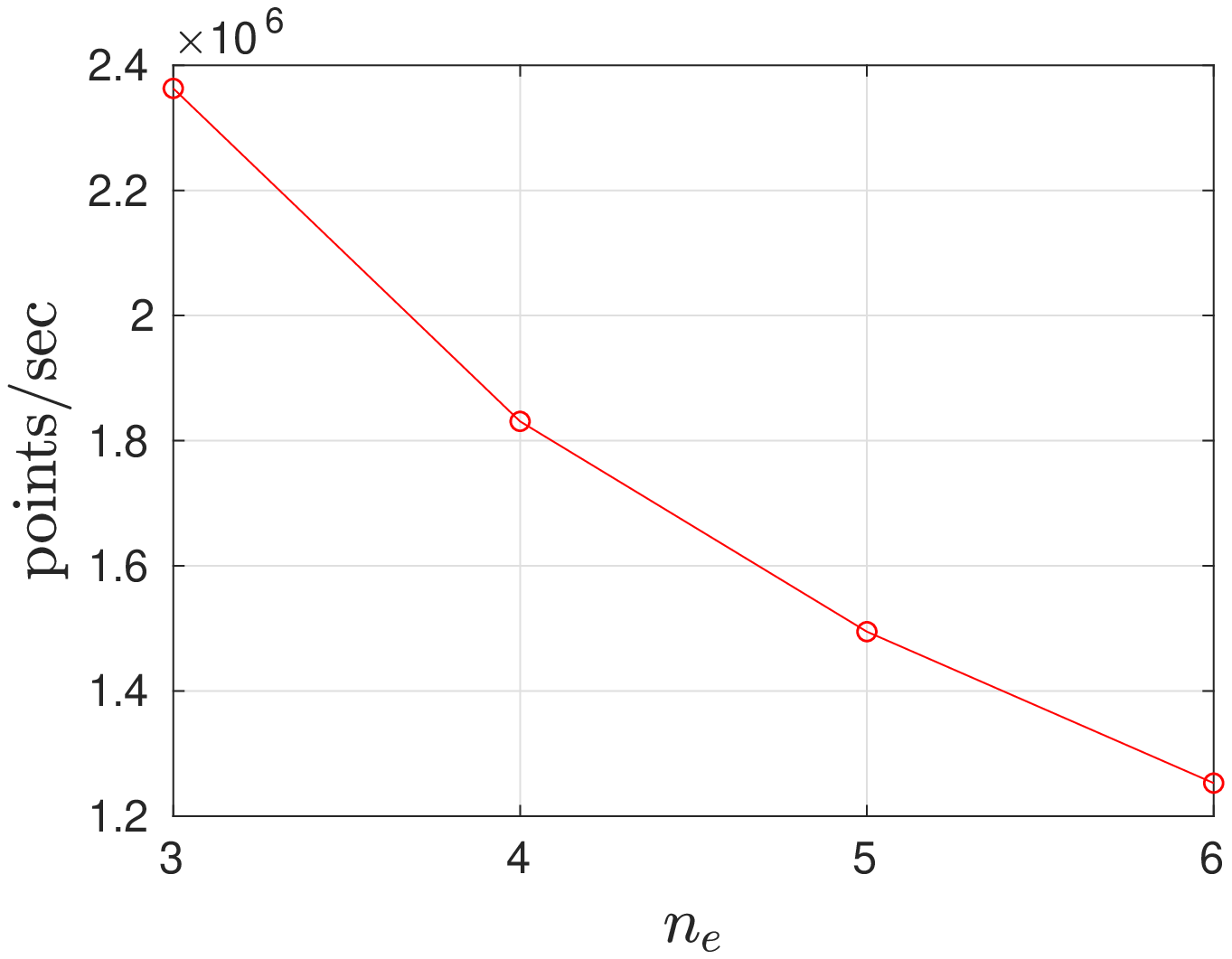}
\caption{Similar to \cref{fig7}, but targets and sources are different.
  Left: Throughput of the FGT without any precomputation.
  Right: Throughput of the FGT when points are pre-sorted.
}
\label{fig9}
\end{figure}

Next we discuss the 1D FGT for the case when the target points are different
from source points. The forward pass is described by the {\sc MATLAB} code fragment
in \cref{fig6}. And the backward pass is similar. Besides calculating the
exponentials $e^{-t_k\frac{x_{i+1}-x_i}{\sqrt{\delta}}}$ for march in targets, we also
need to calculate the exponentials with $x_{i+1}-x_i$ replaced by $x_i-y_j$
or $y_j-x_i$. We tested the code for uniform distribution on $[0,1]$
and the results are shown in \cref{table2} and \cref{fig9}. 
In \cref{table2}, $t_{\rm linear}$ is the computational time excluding
the sorting step. As compared to \cref{table1}, the sorting time
is doubled since both target and source points need to be sorted.
$t_{\rm linear}$ is almost doubled since there are twice as many number
of exponentials needed to be calculated. \Cref{fig9} shows the throughput
of the FGT, where the data points are obtained similarly as in \cref{fig7}.
We observe that the throughput of the full algorithm ranges from $1.5$ to $0.95$
million points per second for $n_e=3,\ldots,6$. We remark that the throughput
of the algorithm will be close to the right panel of \cref{fig7} if we precompute
and store all needed exponentials.
\begin{table}[!ht]
  \caption{Results of the 1D FGT on different target and source points.
  }
\setlength\tabcolsep{3pt}
  \sisetup{
  tight-spacing=true
  }
\centering
\begin{tabular}{S[group-separator = {,},table-format=8.0e0]
    S[table-format=3.0e0]
    S[scientific-notation = fixed,fixed-exponent = 0]
    S[scientific-notation = fixed,fixed-exponent = 0]
    S[scientific-notation = fixed,fixed-exponent = 0]
    S[scientific-notation = true,table-format=3.1e-2]}
  \toprule
  ${N=M}$  & ${\,\,\,\,\,\,\,n_{\rm e}}$ & ${t_{\rm sort}}$ & ${\,\,\,\,\,\,\,t_{\rm linear}}$ & ${t_{\rm total}}$ & {Error}  \\
  \midrule
    100000 & 3 & 0.160D-01 & 0.440D-01 & 0.640D-01 & 0.44D-05 \\
   1000000 & 3 & 0.188D+00 & 0.432D+00 & 0.664D+00 & 0.44D-05 \\
  10000000 & 3 & 0.205D+01 & 0.418D+01 & 0.675D+01 & 0.43D-05 \\
&&&&&\\
    100000 & 4 & 0.120D-01 & 0.520D-01 & 0.680D-01 & 0.56D-07 \\
   1000000 & 4 & 0.188D+00 & 0.556D+00 & 0.856D+00 & 0.55D-07 \\
  10000000 & 4 & 0.204D+01 & 0.541D+01 & 0.800D+01 & 0.55D-07 \\
&&&&&\\
    100000 & 5 & 0.160D-01 & 0.680D-01 & 0.760D-01 & 0.42D-08 \\
   1000000 & 5 & 0.172D+00 & 0.672D+00 & 0.876D+00 & 0.62D-09 \\
  10000000 & 5 & 0.204D+01 & 0.676D+01 & 0.918D+01 & 0.57D-09 \\
&&&&&\\
    100000 & 6 & 0.120D-01 & 0.760D-01 & 0.920D-01 & 0.79D-11 \\
   1000000 & 6 & 0.172D+00 & 0.792D+00 & 0.101D+01 & 0.68D-11 \\
  10000000 & 6 & 0.220D+01 & 0.831D+01 & 0.110D+02 & 0.10D-09 \\

 \bottomrule
 \end{tabular}
\label{table2}
\end{table}

\section{Conclusions and further discussions}\label{sec:conclusion}
Similar to other functions discussed in \cite{schmelzer2007etna,trefethen2006bit},
best rational approximations to the exponential
function on $\mathbb{R}^-$ may be used to construct nearly optimal
SOE approximations for the Gaussian kernel, 
and a simple fast Gauss transform in one dimension has been built
upon such approximations.
The resulting algorithm is benign to parallelization. Indeed, the
most expensive step of the algorithm is computing all needed exponentials,
which is trivial to parallelize. And forward and backward passes
can be easily parallelized over different modes.

The SOE approximations may also be used to speed up the FGTs in higher
dimensions. Instead of standard Fourier spectral expansion, the SOE
approximations can be used to diagonalize translation operators.
The algorithm will be similar to the one in \cite{cheng1999jcp}.

Another interesting fact is about the ill-conditioning of forward and
backward passes. The numerical experiments in \Cref{sec:experiments}
show that the ill-conditioning of the inverse Laplace transform can be overcome
and highly accurate (to almost machine precision) and efficient SOE approximations
can be obtained. However, when we use such SOE approximation for the FGT,
the forward and backward passes may further reduce the accuracy of the overall
calculation. Indeed, the 15-digit accurate SOE approximation discussed
in \cref{sec:roundoff} gives about 11-12 digit accurate solution of the Gauss
transform when $N$ is one million. The SOE approximations obtained by best
rational approximations do not suffer from any ill-conditioning for $n$ up to $12$,
or equivalently, about ten-digit accuracy.
It would be nice if one could develop 
accurate and efficient SOE approximations without suffering of the ill-conditioning
of forward and backward passes beyond that regime.

\section*{Acknowledgments}

The author was supported by the National Science Foundation under grant
DMS-1720405, and by the Flatiron Institute, a division of the Simons
Foundation.  

\bibliographystyle{abbrv}
\bibliography{journalnames,fgt}

\end{document}